\documentclass[11pt]{amsart}
\usepackage{amsmath}
\usepackage{a4wide}
\usepackage[utf8]{inputenc}
\usepackage{amssymb}
\usepackage{amsopn}
\usepackage{epsfig}
\usepackage{amsfonts}
\usepackage{latexsym}
\usepackage{amsthm}
\usepackage{enumerate}
\usepackage[UKenglish]{babel}
\usepackage{bm}
\usepackage{verbatim}
\usepackage{color}
\usepackage{pgf}
\usepackage{tikz}

\usepackage{mathrsfs}   

\usetikzlibrary{arrows,automata}

\newtheorem{theorem}{Theorem}[section]
\newtheorem{lemma}[theorem]{Lemma}
\newtheorem{proposition}[theorem]{Proposition}

\theoremstyle{definition}
\newtheorem{definition}[theorem]{Definition}

\theoremstyle{remark}
\newtheorem{remark}[theorem]{Remark}

\numberwithin{equation}{section}

\newcommand{\R}{\ensuremath{\mathbb{R}}}
\newcommand{\Q}{\ensuremath{\mathbb{Q}}}
\newcommand{\N}{\ensuremath{\mathbb{N}}}
\newcommand{\nn}{\ensuremath{\mathcal{N}}}

\renewcommand{\c}{ {\mathbf{c}}}
\renewcommand{\d}{ {\mathbf{d}}}
\renewcommand{\i}{ {\mathbf{i}}}
\renewcommand{\j}{ {\mathbf{j}}}

\renewcommand{\u}{\ensuremath{\mathbf{u}}}

\renewcommand{\v}{\mathbf{v}}

  \newcommand{\flr}[1]{\left\lfloor#1\right\rfloor}

\newcommand{\set}[1]{\left\{#1\right\}}
\newcommand{\la}{\lambda}
\newcommand{\ga}{\gamma}

\newcommand{\ep}{\varepsilon}
\newcommand{\f}{\infty}

\newcommand{\Om}{\Omega}

\newcommand{\si}{\sigma}

\begin{document}

\title[Projections of four corner Cantor set]{Projections of four corner Cantor set: total self-similarity, spectrum and unique codings}

\author[D. Kong]{Derong Kong}
\address[D. Kong]{College of Mathematics and Statistics, Center of Mathematics, Chongqing University, Chongqing, 401331, P.R.China.}
\email{derongkong@126.com}

\author[B. Sun]{Beibei Sun}
\address[B. Sun]{College of Mathematics and Statistics,   Chongqing University, Chongqing, 401331, P.R.China.}
\email{beibeisun9@163.com}

\date{\today}
\dedicatory{}



\begin{abstract}
Given $\rho\in (0,1/4]$, the four corner Cantor set $E\subset \mathbb{R}^{2}$ is a self-similar set generated by the iterated function system
\[
  \left\{(\rho x, \rho y), \quad(\rho x, \rho y+1-\rho),\quad (\rho x+1-\rho, \rho y),\quad(\rho x+1-\rho,\rho y+1-\rho)\right\}.
  \]
  For $\theta\in[0,\pi)$ let $E_\theta$ be the  orthogonal projection of $E$ onto a line with an angle $\theta$ to the $x$-axis.
  In this paper we give a complete characterization on which the projection $E_\theta$  is totally self-similar.
  We also study the spectrum of $E_\theta $, which turns out that the spectrum of $E_\theta$ achieves its maximum value if and only if $E_\theta $ is totally self-similar.
Furthermore, when $E_\theta$ is totally self-similar, we calculate its Hausdorff dimension and study the subset $U_\theta $ which consists of all $x\in E_\theta $ having a unique coding.
In particular, we show that $\dim_H U_\theta=\dim_H E_\theta$ for Lebesgue  almost every $\theta \in[0,\pi)$.
Finally, for $\rho=1/4$ we describe the distribution of $\theta $ in which $E_\theta$ contains an interval. It turns out that the possibility for $E_\theta$ to contain an interval is smaller than that for $E_\theta$ to have an exact overlap.
\end{abstract}

\maketitle

\section{Introduction}

The study of linear projections of {a planar set} has a long history, {which can be dated}  back to Besicovitch \cite{Besicovitch-1938} and Marstrand \cite{Marstrand-1954}:
for a Borel or analytic set $E\subset\R^2$, let $E_\theta=\mathrm{proj}_\theta(E)$ denote its orthogonal projection {of $E$} onto a line at an angle $\theta$ to the $x$-axis. {Then for Lebesgue almost every $\theta\in[0,\pi)$ we have
$\dim_H E_\theta=\min\set{\dim_H E, 1}$, and in particular, if $\dim_H E>1$ then  $Leb(E_\theta)>0$.
Here $\dim_H$ denotes the Hausdorff dimension. In this paper we study projections of the four corner Cantor set $E$ (cf.~\cite[Ch.~10]{Mattila-2015}), and give a complete characterization for which   $E_\theta$ is totally self-similar (see Definition \ref{def:totoally-self-similar}). Moreover, we study the spectrum of $E_\theta$ (see Definition \ref{def of spec}) and show that $E_\theta$ is totally self-similar if and only if its spectrum achieves its maximum value. Assuming $E_\theta$ is totally self-similar, we calculate its Hausdorff dimension and study its subset $U_\theta$ which consists of all $x\in E_\theta$ having a unique coding. We show that $\dim_H U_\theta=\dim_H E_\theta$ for Lebesgue almost every $\theta\in[0,\pi)$. Furthermore, when $\rho=1/4$ we give the distribution of $\theta$ in which $E_\theta$ contains an interval.
}

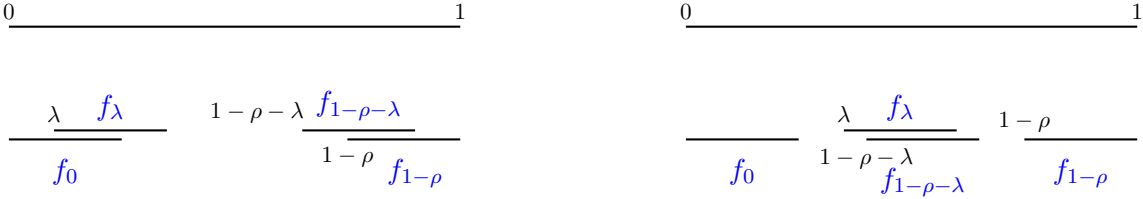
\begin{figure}[h!]
\begin{center}
\begin{tikzpicture}[
scale=6,
axis/.style={very thick, ->},
important line/.style={thick},
dashed line/.style={dashed, thin},
pile/.style={thick, ->, >=stealth', shorten <=2pt, shorten
>=2pt},
every node/.style={color=black}
]
\draw[important line] (0, -0.05)--(1, -0.05);
\node[above,scale=0.8pt] at(0, -0.05){$0$};\node[above,scale=0.8pt] at(1, -0.05){$1$};
 \draw[important line] (1.5, -0.05)--(2.5, -0.05);
 \node[above,scale=0.8pt] at(1.5, -0.05){$0$};\node[above,scale=0.8pt] at(2.5, -0.05){$1$};

\draw[important line] (0, -0.3)--(0.25, -0.3); \draw[important line] (0.1, -0.28)--(0.35, -0.28);
\node[above,scale=0.8pt] at(0.1, -0.28){$\lambda $};\node[blue,below,scale=1pt] at(0.125, -0.32){$f_0$};
\node[blue,above,scale=1pt] at(0.225, -0.28){$f_{\lambda }$};

 \draw[important line] (0.75, -0.3)--(1, -0.3); \draw[important line] (0.65, -0.28)--(0.9, -0.28);
\node[above,scale=0.8pt] at(0.55, -0.28){$1-\rho-\lambda$};\node[below,scale=0.8pt] at(0.75, -0.3){$1-\rho$};
\node[blue,above,scale=1pt] at(0.775, -0.28){$f_{1-\rho-\lambda}$};\node[blue,below,scale=1pt] at(0.9, -0.32){$f_{1-\rho}$};

 \draw[important line] (1.5, -0.3)--(1.75, -0.3); \draw[important line] (1.85, -0.28)--(2.1, -0.28);
\node[above,scale=0.8pt] at(1.85, -0.28){$\lambda $};\node[blue,below,scale=1pt] at(1.625, -0.32){$f_0$};
\node[blue,above,scale=1pt] at(1.975, -0.28){$f_{\lambda }$};

\draw[important line] (2.25, -0.3)--(2.5, -0.3); \draw[important line] (1.9, -0.3)--(2.15, -0.3);
\node[below,scale=0.8pt] at(1.9, -0.3){$1-\rho-\lambda$};\node[above,scale=0.8pt] at(2.25, -0.3){$1-\rho$};
\node[blue,below,scale=1pt] at(2.025, -0.34){$f_{1-\rho-\lambda}$};\node[blue,below,scale=1pt] at(2.375, -0.32){$f_{1-\rho}$};
\end{tikzpicture}
\end{center}
\caption{The first two levels for the geometric construction of $E_\lambda$ with $\lambda\in (0,\rho)$ (left) and $\lambda\in \left(\frac{1-2\rho}{2},\frac{1-\rho}{2}\right)$ (right).}\label{Fig:1}
\end{figure}

Given $\rho\in(0,1/4]$,
  let $E\subset\R^2$ be the four corner Cantor set, which is a self-similar set   generated by the \emph{iterated function system} (IFS)
 \[
\set{(\rho x, \rho y), \quad(\rho x, \rho y+1-\rho),\quad (\rho x+1-\rho, \rho y),\quad(\rho x+1-\rho,\rho y+1-\rho)}.
\]
It is well known that $\dim_H E=\frac{2\log 2}{-\log \rho}$ and its Hausdorff measure $\mathcal H^{\frac{2\log 2}{-\log \rho}}(E)\in(0,\f)$ (cf.~\cite{Hutchinson_1981}).
Furthermore, the self-similar set $E$ can be written algebraically as

\begin{equation}\label{IFSofE}
E=\set{\sum_{i=1}^\infty \rho^{i-1} d_i: d_i\in\set{(0,0), (0,1-\rho),(1-\rho,0),(1-\rho,1-\rho)}~\forall i\in \mathbb{N}}.
\end{equation}
For $t\in\mathbb{R}$ let $E(t)$ be its orthogonal projection onto a line with slope $t$.  Then
\[E(t)=\set{\frac{x+t y}{\sqrt{1+t^2}}: (x,y)\in E},\]
 and it is also a self-similar set.
By symmetry and scaling, we may reduce the projection $E(t)$   to    the self-similar set $E_\la$ generated by the IFS
\begin{equation}\label{eq:IFS}
\mathcal F_\lambda :=\set{f_{d}(x)=\rho x+d: d\in \Omega_{\lambda}}\quad \textrm{with}\quad\Omega_{\lambda}:=\set{0,\lambda ,1-\rho-\lambda,1-\rho},
\end{equation}  where $\la\in[0,1-\rho]$.
In other word, it suffices to consider
\begin{equation}\label{di}
E_\lambda=\bigcup_{d\in\Omega_\lambda} f_d(E_\lambda)=\left\{\sum_{i=1}^{\infty} \rho ^{i-1}d_i  \mid  d_i\in \Omega_\lambda~ \forall i\in \mathbb{N} \right\},\quad \la\in[0,1-\rho].
\end{equation}
{Note that  $\rho\in(0,1/4]$ is fixed, and   it is always suppressed  in our notation.}
In fact, we can   restrict our parameter $\lambda $ to the interval $(0,\rho)\cup \left(\frac{1-2\rho}{2},\frac{1-\rho}{2}\right)$.
Note  by the symmetry that $E_\lambda$  has the same geometrical structure as $E_{1-\rho-\lambda}$ for any $\lambda\in [0,1-\rho)$. Then we only need to consider $\lambda\in [0,\frac{1-\rho}{2}]$.
Moreover, for $\lambda\in [\rho,\frac{1-2\rho}{2}]$ the self-similar set $E_\lambda$ satisfies the open set condition, and for $\lambda=0$ or $\lambda=\frac{1-\rho}{2}$ the self-similar set $E_\lambda$ can be degenerated to a self-similar set satisfies the {\emph{strong separation condition} (SSC)}.
So we only need  to consider  $\lambda \in(0,\rho)\cup \left(\frac{1-2\rho}{2},\frac{1-\rho}{2}\right)$, and in this case the self-similar set $E_\lambda$ has non-trivial overlaps (see Figure \ref{Fig:1} for the two types of overlapping structure).

In 2004 Broomhead, Montaldi and Sidorov \cite{Bro-Mon-Sid-04} introduced the following finer family of self-similar sets with overlaps. For $n\in \mathbb{N}\cup \left\{0\right\}$ let $\Omega_\lambda^n:=\left\{i_1i_2\cdots i_n\mid i_k\in \Omega_\lambda, 1\leq k\leq n \right\}$,
where for $n=0$ we set $\Omega_\lambda^0:=\{\epsilon\}$ with $\epsilon$ the empty word.
Let  $\Omega_\lambda^*$ be the set of all finite words over the alphabet $\Omega_\lambda$, i.e., $\Omega_\lambda^*=\bigcup_{n=0}^{\infty} \Omega_\lambda^n$.
Furthermore, let $\Omega_\lambda^{\mathbb{N}}$ be the set of all infinite sequences over the alphabet $\Omega_\lambda$.
For $\i=i_1i_2\cdots i_n\in \Omega_\lambda^*$ we write
\begin{equation}\label{eq:composition-maps}
f_{\i}(x):=f_{i_1} \circ f_{i_2}\circ \cdots  \circ f_{i_n}(x)=\rho^n x+\sum_{k=1}^n\rho^{k-1}i_k
\end{equation} as compositions of maps.
In particular, for $\i=\epsilon$ we set $f_{\epsilon}$ as the identity map.

\begin{definition}\label{def:totoally-self-similar}
Let $I=[0,1]$ be the convex hull of $E_\la$. The  self-similar set
$E_\lambda$ is called \emph{totally self-similar} if
$$
f_{\mathbf{i}}\left(E_\lambda\right)=f_{\mathbf{i}}(I) \cap E_\lambda \quad \text { for any } \mathbf{i} \in \Omega_\lambda^*.
$$
\end{definition}
Our first result characterizes when $E_\lambda$ is totally self-similar.
For $k\in \mathbb{N}\cup\left\{0\right\}$ we define
\begin{equation}\label{eq:lambda-gamma-eta}
  \lambda_k:=\frac{ \rho(1-\rho ^k)}{1+\rho ^k},\quad \gamma_k:=\rho(1-\rho^k)\quad\textrm{and} \quad \eta_k:=\frac{1-2\rho+\rho^{k+1}}{2}.
\end{equation}
Then it is clear that
\begin{align*}
&0=\lambda_0=\gamma_0<\lambda_1<\gamma_1<\cdots<\lambda_k<\gamma_k<\cdots<\rho,  \\
&\frac{1-\rho}{2}=\eta_0>\eta_1>\eta_2>\cdots>\eta_k>\eta_{k+1}>\cdots>\frac{1-2\rho}{2}.
\end{align*}
 {Furthermore, $ \lambda_k,\gamma_k\nearrow \rho$ and $\eta_k \searrow \frac{1-2\rho}{2}$ as $k\rightarrow \infty$.}
\begin{theorem}\label{th:totally-self-similar}
Let $\lambda \in(0,\rho)\cup \left(\frac{1-2\rho}{2},\frac{1-\rho}{2}\right)$.
\begin{enumerate}
[{\rm(i)}]
\item If $\lambda \in(0,\rho)$, then $E_\lambda$ is totally self-similar if and only if $\lambda=\lambda_k$ or $\gamma_k$ for some $k\in \mathbb{N}$.
\item If $\lambda \in\left(\frac{1-2\rho}{2},\frac{1-\rho}{2}\right)$, then $E_\lambda$ is totally self-similar if and only if $\lambda=\eta_k$ for some $k\in \mathbb{N}$.
\end{enumerate}

\end{theorem}

Given  two   words $\mathbf{i}, \mathbf{j}\in\Omega_\la^n$, by (\ref{eq:composition-maps})  it is clear that $f_{\mathbf{i}}=f_{\mathbf{j}}$ if and only if $f_{\mathbf{i}}(0)=f_{\mathbf{j}}(0)$.
So, the scaled distance $\frac{\left|f_{\mathbf{i}}(0)-f_{\mathbf{j}}(0)\right|}{\rho ^n}$ describes the closeness of the two maps $f_{\mathbf{i}}$ and $f_{\mathbf{j}}$, which reveals  the complexity of the overlapping structure of $E_\lambda$. Let
\begin{align*}
A_\lambda&:= \left\{\frac{\left|f_{\mathbf{i}}(0)-f_{\mathbf{j}}(0)\right|}{\rho ^{n}}: \mathbf{i}, \mathbf{j} \in \Omega_\lambda^n \text { with } f_{\mathbf{i}} \neq f_{\mathbf{j}};~ n\in\N  \right\}\\
&=\left\{\left|\sum_{i=1}^{n} \frac{d_i}{\rho^{i}}\right| \neq 0: d_i \in\Omega_\lambda^{\pm };~ n\in\N  \right\},
\end{align*}
where {$\Omega_\lambda^{\pm }:=\Omega_\lambda-\Omega_\lambda=\left\{0,\pm \lambda,\pm(1-\rho-2\lambda),\pm(1-\rho-\lambda),\pm(1-\rho)\right\}$.}

Motivated by the spectrum from non-integer base expansions (cf.~\cite{Erdos_Joo_Komornik_1990}) and the spectrum for overlapping self-similar sets (cf.~\cite{Dajani-Kong-Yao-2019}) we consider the spectrum of $E_\la$.

\begin{definition}\label{def of spec}
For $\lambda \in(0,\rho)\cup \left(\frac{1-2\rho}{2},\frac{1-\rho}{2}\right)$, the \emph{spectrum} of $E_\lambda$ is defined by
$$
l_\lambda:=\inf  A_\lambda=\inf \left\{\left|\sum_{i=1}^{n} \frac{d_i}{\rho^{i}}\right| \neq 0: d_i \in\Omega_\lambda^{\pm };~ n\in \mathbb{N}\right\}.
$$
\end{definition}

Our second result describes the spectrum $l_\lambda$ of $E_\la$.

\begin{theorem}\label{spectrum}\mbox{}
\begin{enumerate}[\rm(i)]
\item For any $k\in\N$ we have
\[
l_{\lambda_k}=1-\rho-\lambda_k,\quad l_{\gamma_k}=1-\rho\quad\textrm{and}\quad l_{\eta_k}=1-\rho.
\]
  \item If $\lambda \in(0,\rho)$, then $E_\lambda$ is not totally self-similar if and only if $l_{\lambda}< 1-\rho-\lambda$.
 \item If $\lambda \in\left(\frac{1-2\rho}{2},\frac{1-\rho}{2}\right)$, then $E_\lambda$ is not totally self-similar if and only if $l_{\lambda}< 1-\rho$.
\end{enumerate}
\end{theorem}

\begin{remark}
\begin{enumerate}[{\rm(i)}]
  \item By Theorem \ref{spectrum} it follows that the spectrum $l_\la$ attains its maximum value if and only if $E_\la$ is totally self-similar.

  \item
  If $0<\rho<1/9$, then by (\ref{IFSofE}) we have $\dim_H (E-E)<1$.
  By \cite[Lemma 2.7]{Peres-2000} it follows that $\mathcal H^s(E_\la)>0$ for Lebesgue  almost every $\lambda\in (0,\rho)\cup (\frac{1-2\rho}{2},\frac{1-\rho}{2})$, where $s=\frac{\log 4}{-\log\rho}$.
 Hence, by \cite[Corollary 3.2]{Farka-2015} we can deduce that $E_\lambda$ satisfies the weak separate condition for Lebesgue almost every $\lambda$. Note that $l_\lambda>0$ is equivalent to that $E_\la$ satisfies the weak separation condition (cf.~\cite{zerner-1996, Farka-2015}). So, if $\rho\in(0, 1/9)$ then   $l_\lambda>0$ for Lebesgue almost every $\la$.
 \end{enumerate}
\end{remark}

In view of (\ref{di}) and (\ref{eq:composition-maps}), for each $x \in E_\lambda$ we can find a sequence $\left(d_i\right)=d_1d_2\ldots\in \Omega_{\lambda}^{\mathbb{N}}$ such that
\begin{equation}\label{eq:coding-map}
x=\lim _{n \rightarrow \infty} f_{d_1 \cdots d_n}(0)=\sum_{i=1}^{\infty} \rho ^{i-1}d_i=:\pi_\la((d_i)).
\end{equation}
The infinite sequence $\left(d_i\right)$ is called a \emph{coding} of $x$ with respect to the digit set $\Omega_\lambda$.
Since $E_\la$ has overlaps,  $x\in E_\la$ might have multiple codings.
In this paper, we are also interested in the subset
\[
U_\la:=\set{x\in E_\la: \#\pi_\la^{-1}(x)=1}.
\]
Then each $x\in U_\la$ has a unique coding.

 Our third result shows that if $E_\la$ is  {totally self-similar}, then $\dim_H U_\la<\dim_H E_\la$. Furthermore, we give the analytic formula for the dimension of $U_\la$.
\begin{theorem}\label{th:unique-coding-totally-self-similar}
If $E_\la$ is totally self-similar, i.e., $\la\in\bigcup_{k=1}^\f\set{\la_k,\ga_k,\eta_k}$, then
\[
\dim_H U_\la<\dim_H E_\la.
\]
\begin{enumerate}
  [{\rm(i)}]
  \item If $\la=\la_k$ for some $k\in\N$, then $\dim_H U_\la=s$, where $s\in(0,1)$ is an appropriate root of
  \[4\rho^{s}-2\rho^{ks}=1.\]

  \item If $\la=\ga_k$ for some $k\in\N$, then $\dim_H U_\la=s, \dim_H E_\la=t$, where $s, t\in(0,1)$ are respectively appropriate roots of
  \[4\rho^s-\rho^{ks}=1\quad\textrm{and}\quad 4\rho^t-2\rho^{(k+1)t}=1.\]

  \item If $\la=\eta_k$ for some $k\in\N$, then $\dim_H U_\la=s, \dim_H E_\la=t$, where $s, t\in(0,1)$ are respectively  appropriate roots of
  \[4\rho^s-2\rho^{(k+1)s}=1\quad\textrm{and}\quad 4\rho^t-\rho^{(k+1)t}=1.\]
\end{enumerate}
\end{theorem}
\begin{remark}
  When $\la=\la_k$, although the self-similar set $E_\la$ can be represented as a graph-directed set,   the directed graph does not satisfy  the open set condition. So we don't know how to calculate  the Hausdorff dimension of $E_\la$ in this case.
\end{remark}
If $E_\la$ is totally self-similar, then   $E_\la$ has exact overlaps, i.e., $f_{\i}=f_{\j}$ for some $\i\ne\j$ (see Remark \ref{rem:41} and Lemma \ref{lem:f0-cap-flambda}). This implies that $\dim_H E_\la<\dim_S E_\la=\frac{\log 4}{-\log\rho}$, where $\dim_S$ denotes the similarity dimension. However,   by  \cite[Theorem 2.1]{Shmerkin-2015} we know that $\dim_H E_\la=\frac{\log 4}{-\log\rho}$ for Lebesgue almost every $\lambda \in(0,\rho)\cup\left(\frac{1-2\rho}{2},\frac{1-\rho}{2}\right)$.
Our fourth result  states that for typical $\lambda$ the  {univoque} set $U_\la$ has the same Hausdorff dimension as $E_\la$.
\begin{theorem}
  \label{th:unique-codings-typical}\mbox{}
  If $\rho\in(0,1/4)$, then for Lebesgue almost every $\lambda \in(0,\rho)\cup\left(\frac{1-2\rho}{2},\frac{1-\rho}{2}\right)$ we have
  \[
  \dim_H U_\la=\dim_H E_\la=\frac{ \log 4}{-\log \rho}.
  \]
In particular, if $\rho\in(0,1/16)$, then
  $
  U_\la=E_\la
  $ for Lebesgue almost every $\la\in(0,\rho)\cup\left(\frac{1-2\rho}{2},\frac{1-\rho}{2}\right)$.
 \end{theorem}
\begin{remark}
  \begin{enumerate}[{\rm(i)}]
 \item Theorem \ref{th:unique-codings-typical} can be easily extended to all $\la\in\R$;
 \item If $\rho\in(0,1/16)$, then Theorem \ref{th:unique-codings-typical} suggests that for Lebesgue almost every $\la\in(0,\rho)\cup\left(\frac{1-2\rho}{2},\frac{1-\rho}{2}\right)$ the self-similar set $E_\la$ satisfies the {SSC}, i.e., $f_i(E_\la)\cap f_j(E_\la)=\emptyset$ for any $i\ne j\in\Omega_\lambda$. {An extension to a larger class of self-similar sets with overlaps  can be found in \cite{Baker-Kong-Wang-2024}.}
  \end{enumerate}
\end{remark}

{When $\rho=1/4$,   projections of the four corner Cantor set $E$ defined in (\ref{IFSofE}) are extensively studied  (cf.~\cite{Kenyon-1997, Rao-Wen-1998, Peres-2000, Mattila-2015}).   Note by (\ref{eq:IFS}) that   the scaled projection $E_\la$ is a self-similar set generated by the IFS
\[
\set{\hat f_d(x)=\frac{x+d}{4}: d\in\set{0,4\la, 3-4\la, 3}}.
\]
It is known that $E_\la$ has zero Lebesgue measure for Lebesgue almost every $\la\in\R$. Indeed, there are only countably many $\la$ for which $E_\la$ has positive Lebesgue measure. More precisely, if $\la\notin\mathbb Q$, then $E_\la$ has zero Lebesgue measure but full Hausdorff dimension; if $\la\in\mathbb Q$ and $E_\la$ has   exact overlaps (see its definition below), then $\dim_H E_\la<1$;  if $\la\in\mathbb Q$ but $E_\la$ does not have an exact overlap, then $E_\la$ is a perfect set containing a non-degenerate interval.
}
\begin{definition}
  \label{def:exact-overlap}
  $E_\la$ is said to have an \emph{exact overlap} if there exist two blocks $\i=i_1\ldots i_n, \j=j_1\ldots j_n\in\set{0,4\la, 3-4\la, 3}^n$ such that
  $
  \hat f_{\i}=\hat f_{\j}.
  $
\end{definition}
The following  complete characterization of exact overlaps  of $E_\la$ can be essentially deduced from \cite[Theorem 10.5]{Mattila-2015} (see also, \cite{Kenyon-1997, Rao-Wen-1998, Peres-2000}). For   $n\in\N$ let $ord_2(n)$ be the highest power of $2$ that divides $n$ (cf.~\cite[P.~2]{Koblitz-1984}).
\begin{theorem} [\cite{Mattila-2015}]
  \label{th:martilla}
 Let $\rho=1/4$ and $\la\in(0,\frac{3}{8})$. Then $E_\la$ has an exact overlap if and only if
 \[\la=\frac{3p}{4(p+q)}\in\Q\quad\textrm{ with }\quad (p,q)\in W,
 \] where
\begin{equation}\label{eq:W}
 W:=\set{(p,q)\in\N^2:\textrm{both }ord_2(p)\textrm{ and }ord_2(q)\textrm{ are even},~ p<q\textrm{ and }p,q\textrm{ are coprime}}.
 \end{equation}
 Furthermore, the following statements hold true.
 \begin{enumerate}
   [{\rm(i)}]
   \item If $\la=\frac{3p}{4(p+q)}\in\Q$ in reduced form with $(p,q)\in W$, then $\dim_H E_\la<1$.
   \item If $\la=\frac{3p}{4(p+q)}\in\Q$ in reduced form with $(p,q)\notin W$, then $E_\la$ contains an interval.
   \item If $\la\notin\Q$, then $C_\la$ has zero Lebesgue measure and $\dim_H E_\la=1$.
 \end{enumerate}
\end{theorem}
\begin{remark}
  Note that in \cite[Theorem 10.5]{Mattila-2015} the result was stated using the following notation. For $n\in\N$ let $n^*\in\set{1,2,3}$ be defined by
  \[
  n^*=\frac{n}{4^{j_0}}\mod 4,
  \]
  where $j_0$ is the largest integer $j$ such that $4^j$ divides $n$. One can easily verify that $n^*$ is odd if and only if $ord_2(n)$ is even. So, Theorem \ref{th:martilla} is the same as  \cite[Theorem 10.5]{Mattila-2015}.
\end{remark}



Our final result describes the density of $W$ in $\N^2$, which reveals the possibility in which the projection $E_\la$ has an exact overlap. We also consider the density of
\begin{equation}\label{eq:hat-W}
\hat W:=\set{(p,q)\in\N^2: ord_2(p)\textrm{ odd  or }ord_2(q) \textrm{ odd,  and }p<q\textrm{ with }p, q\textrm{ coprimes}},
\end{equation}
 which describes the possibility in which  $E_\la$ contains a non-degenerate interval. For a set $A$ let $\#A$ denote its cardinality.
\begin{theorem}
  \label{th:density-exact-overlap}
  Let $W$ and $\hat W$ be defined as in (\ref{eq:W}) and (\ref{eq:hat-W}) respectively.  Then
  \[
  \lim_{N\to\f}\frac{\#(W\cap[1, N]^2)}{N^2}=\frac{5}{3\pi^2}\quad\textrm{and}\quad
  \lim_{N\to\f}\frac{\#(\hat W\cap[1, N]^2)}{N^2}=\frac{4}{3\pi^2}.
  \]
\end{theorem}

Theorem \ref{th:density-exact-overlap} indicates that the possibility for $E_\la$ to contain a non-degenerate interval is  smaller than that for $E_\la$ to have an exact overlap.

The rest of the paper is organized as follows. In the next section we give a complete characterization when $E_\la$ is totally self-similar and prove Theorem \ref{th:totally-self-similar}. In Section \ref{sec:spectrum} we study the spectrum of $E_\la$ and prove Theorem \ref{spectrum}. In particular, we show that $E_\la$ is totally self-similar if and only if its spectrum achieves its maximum value. In Section \ref{sec:unique-codings-totally-self-similar} we consider the subset $U_\la$ which consists of all $x\in E_\la$ having a unique coding, and calculate its Hausdorff dimension (Theorem \ref{th:unique-coding-totally-self-similar}). In Section \ref{sec:typical-dimension-U} we show that $\dim_H U_\la=\dim_H E_\la$ for typical $\la$, and prove Theorem \ref{th:unique-codings-typical}. Finally, we consider the four corner Cantor set $E$ with dimension one, i.e., $\rho=1/4$. Although the projection  $E_\la$ was extensively studied in the literature (see, e.g., \cite[Ch.~10]{Mattila-2015}), we add a new result on the distribution of $\la$ in which $E_\la$ contains an interval.

\section{Total self-similarity of $E_\lambda$}\label{sec:totoally-self-similarity}
In this section we will characterize the total self-similarity of $E_\lambda$, and prove Theorem \ref{th:totally-self-similar}.
Given $\rho\in(0,1/4]$ and $\la\in(0,\rho)\cup(\frac{1-2\rho}{2}, \frac{1-\rho}{2})$,  we recall from Definition \ref{def:totoally-self-similar} that $I =[0, 1]$ is the convex hull of $E_\lambda$. Set $I _0 = I$, and for $n\geq 1$ let
$$
I _n:=\bigcup _{\mathbf i \in \Omega_\lambda^n}{f_ \mathbf i(I )},
$$
where $\Om_\la^n$ consists of all length $n$ words over   $\Om_\la=\set{0,\la, 1-\rho-\la,1-\rho}$.
Let
\[H:=I\setminus I_1=I\setminus\bigcup_{d\in\Om_\la}f_d(I)\] be  {a} hole of $E_\lambda$.
The following characterization of total  {self-similarity} of $E_\lambda$ can be found in \cite[Proposition 2.1]{Dajani-Kong-Yao-2019}.
\begin{proposition}\label{equ state}
  The set $E_\lambda$ is totally self-similar if and only if for any two  words $\mathbf{i}, \mathbf{j}\in \Omega_{\lambda}^{n}$ with $n\in\N$,
 \[\textrm{either}\quad f_{\mathbf{i}}=f_{\mathbf{j}}\quad\textrm{or}\quad f_{\mathbf{i}}(I) \cap f_{\mathbf{j}}(H)=\emptyset.\]
\end{proposition}

Recall from (\ref{eq:lambda-gamma-eta}) the definitions of $\la_k, \ga_k$ and $\eta_k$ for $k\in\N$. Let
\begin{equation}\label{ell}
\ell_k :=\frac{\rho(1-\rho^{k+1})}{1+\rho^k+\rho^{k+1}},\quad k\in\N.
\end{equation}
Then by using  {$\rho\in (0,1/4]$} it follows that
\[\lambda_k<\ell_k<\gamma_k\quad \textrm{for all }k\in \mathbb{N}.\]
To prove Theorem \ref{th:totally-self-similar}   we need the following two lemmas.
\begin{lemma}\label{inequ}
Let $\rho\in(0,1/4]$ and $k\in \mathbb{N}$.
  \begin{enumerate}[{\rm(i)}]
 \item   If $\lambda \ge \lambda_k$, then $  f_{\lambda 0^k}(1-\rho-\lambda)>f_{0(1-\rho)^{k-1}(1-\rho-\lambda)}(0).$
 \item If $\lambda \geq \ell_k$, then $  f_{\lambda 0^k}(1-\rho-\lambda)>f_{0(1-\rho)^{k}}(0).$
 \item If $\lambda<\eta_{k-1}$, then $  f_{\lambda (1-\rho)^{k-1}}(0)<f_{(1-\rho-\lambda)0^{k-1}}(0).$
\end{enumerate}
\end{lemma}

\begin{proof}
Since the proofs of the three items are similar, we only prove (i). Note by (\ref{eq:composition-maps}) that
{\begin{equation*}
  \begin{split}f_{\lambda 0^k}(1-\rho-\lambda)&=\lambda+\rho^{k+1}(1-\rho-\lambda),\quad
  f_{0(1-\rho)^{k-1}(1-\rho-\lambda )}(0)=\rho-\rho^{k+1}-\rho^{k}\lambda.
  \end{split}
\end{equation*}}
Then $  f_{\lambda 0^k}(1-\rho-\lambda)>f_{0(1-\rho)^{k-1}(1-\rho-\lambda)}(0)$ is equivalent to
\[\lambda>\frac{\rho-2\rho^{k+1}+\rho^{k+2}}{1+\rho^k-\rho^{k+1}}.\]
  Since $\lambda \ge\lambda_k=\frac{\rho(1-\rho^{k})}{1+\rho^k}$, it suffices to prove
\[\frac{ \rho(1-\rho ^k)}{1+\rho ^k}>\frac{\rho-2\rho^{k+1}+\rho^{k+2}}{1+\rho^k-\rho^{k+1}},\]
which holds by using $0<\rho\le 1/4$.
\end{proof}
For $\lambda\in (0,\rho)$ we recall that
 {$\Omega_{\lambda}^{\pm}= \left\{0,\pm \lambda,\pm(1-\rho-2\lambda),\pm(1-\rho-\lambda),\pm(1-\rho)\right\}.$}
\begin{lemma}\label{lambda1}
Let $\la=\la_1=\frac{\rho(1-\rho)}{1+\rho}$.
  If $d\in \Omega_{\lambda_1}^{\pm}$ and $d\leq 0$, then
  \[d+\frac{\lambda_1}{\rho}\in \Omega_{\lambda_1}^{\pm}.\]
\end{lemma}
\begin{proof}
 Note that
  \[
   \frac{\lambda_1}{\rho}=\frac{1-\rho}{1+\rho}=(1-\rho)\left(1-\frac{\rho}{1+\rho}\right)=1-\rho-\frac{\rho(1-\rho)}{1+\rho}=1-\rho-\lambda_1.
 \]
This implies that
\begin{align*}
  &0+\frac{\lambda_1}{\rho} =1-\rho-\lambda_1\in \Omega_{\lambda_1}^{\pm},  &-\lambda_1+\frac{\lambda_1}{\rho} =1-\rho-2\lambda_1\in \Omega_{\lambda_1}^{\pm},   \\
  &-(1-\rho-2\lambda_1)+\frac{\lambda_1}{\rho}  =\lambda_1\in \Omega_{\lambda_1}^{\pm},   &-(1-\rho-\lambda_1)+\frac{\lambda_1}{\rho}  =0 \in \Omega_{\lambda_1}^{\pm}, \\
 &-(1-\rho)+\frac{\lambda_1}{\rho}  =-\lambda_1 \in \Omega_{\lambda_1}^{\pm} &
\end{align*}
as desired.
\end{proof}
In the following we split our proof of Theorem \ref{th:totally-self-similar} into two subsections for $\la\in(0,\rho)$ and $\la\in(\frac{1-2\rho}{2}, \frac{1-\rho}{2})$, separately.

\subsection{Total self-similarity of $E_\la$ for $\la\in(0,\rho)$}
Since $\lambda\in (0,\rho)$, the hole $H$ is given by (see the left graph of Figure \ref{Fig:1})
 \[H=I \setminus I _1=(\rho+\lambda, 1-\rho-\lambda).\]
First we show that  $\la\in\bigcup_{k=1}^\f\set{\la_k,\ga_k}$ is necessary for the total self-similarity of $E_\la$.
\begin{proposition}\label{NTSS}
  If $\lambda\in (0,\rho)\setminus \bigcup_{k=1}^{\infty} \left\{ \lambda_k,  \gamma_k \right\}$, then $E_\lambda$ is not totally self-similar.
\end{proposition}

\begin{proof}
  Note by (\ref{eq:lambda-gamma-eta}) and (\ref{ell}) that
  $0<\lambda_k<\ell_k<\gamma_k<\lambda_{k+1}$ for any $k\in \mathbb{N}$, and $\lambda_k\nearrow \rho$ as $k\to \infty$.
So, we only need to prove that $E_\lambda$ is not totally self-similar for any
(\uppercase\expandafter{\romannumeral1}) $\lambda\in (0,\lambda_1)$;
(\uppercase\expandafter{\romannumeral2}) $\lambda\in \bigcup_{k=1}^{\infty} (\lambda_k,\ell_k)$; and
(\uppercase\expandafter{\romannumeral3}) $\lambda\in \bigcup_{k=1}^{\infty} [\ell_k,\gamma_k)\cup(\gamma_k,\lambda_{k+1})$. By Proposition \ref{equ state} it suffices to prove in the three different cases that there exist two words $\i, \j\in\Om_\la^n$ for some $n\in\N$ such that
\begin{equation}\label{eq:Not-totally-self-similar}f_{\i}\ne f_{\j}\quad \textrm{and}\quad f_{\i}(I)\cap f_{\j}(H)\ne\emptyset,
\end{equation} where $I=[0,1]$ and {$H=(\rho+\la, 1-\rho-\la)$}.

Case (\uppercase\expandafter{\romannumeral1})\ $\lambda\in (0,\lambda_1)$.
Take $\mathbf{i}=0$, $\mathbf{j}=\lambda\in \Omega_\lambda$. Since
$
  f_{0}(0)=0<\lambda=f_{\lambda}(0),
$
by  {(\ref{eq:composition-maps})} we have $f_{0}\neq f_{\la}$.
Furthermore, note that
\begin{equation}\label{0lambda}
  f_{0}(I)=[0,\rho] \quad \text{ and } \quad f_{\lambda}(H)=(\rho(\rho+\lambda)+\lambda, \rho (1-\rho-\lambda)+\lambda).
\end{equation}
Then by using  {$0<\lambda<\rho\le 1/4$} it is clear  that $0< \rho (1-\rho-\lambda)+\lambda$.
Since  $\lambda<\lambda_1=\frac{\rho(1-\rho)}{1+\rho}$, we have {$\rho>\rho(\rho+\lambda)+\lambda$}.
This together with (\ref{0lambda}) implies that $f_{0}(I)\cap f_{\lambda}(H)\neq \emptyset$, proving (\ref{eq:Not-totally-self-similar}).

Case (\uppercase\expandafter{\romannumeral2})\ $\lambda\in  (\lambda_k,\ell_k)$ for some $k\in \mathbb{N}$.
Take $\mathbf{i}=\lambda 0^k$ and $\mathbf{j}=0(1-\rho)^{k-1}(1-\rho-\lambda )$.
 Since $\lambda>\lambda_k$, by (\ref{eq:composition-maps}) and (\ref{eq:lambda-gamma-eta}) it follows that
 \[f_{\lambda 0^k}(0)=\lambda>\rho-\rho^{k}\la-\rho^{k+1}=f_{0(1-\rho)^{k-1}(1-\rho-\lambda )}(0),\] which implies $f_{\lambda 0^k}\neq f_{0(1-\rho)^{k-1}(1-\rho-\lambda )}$.
 Moreover,
 by (\ref{eq:composition-maps}) it follows that
\begin{equation}\label{lambda0k=01-rhok-1}
  \begin{split}
  &f_{\lambda 0^k}(H) =\left(\lambda +\rho^{k+1}(\rho+\lambda), \lambda +\rho^{k+1}(1-\rho-\lambda)\right), \\
  &f_{0(1-\rho)^{k-1}(1-\rho-\lambda )}(I)=\left[\rho-\rho^{k}\lambda-\rho^{k+1}, \rho-\rho^{k}\lambda\right] .
  \end{split}
\end{equation}
 Since  $\lambda\in  (\lambda_k,\ell_k)$, by (\ref{ell}) we have {$\lambda +\rho^{k+1}(\rho+\lambda)<\rho-\rho^{k}\lambda$}. Furthermore,  by Lemma \ref{inequ} (i) it follows that
 $$
 \rho-\rho^{k}\lambda-\rho^{k+1}=f_{0(1-\rho)^{k-1}(1-\rho-\lambda )}(0)<f_{\lambda 0^k}(1-\rho-\lambda)=\lambda +\rho^{k+1}(1-\rho-\lambda).
 $$
So, by (\ref{lambda0k=01-rhok-1}) it follows that $f_{\lambda 0^k}(H)\cap f_{0(1-\rho)^{k-1}(1-\rho-\lambda )}(I)\neq \emptyset$, establishing (\ref{eq:Not-totally-self-similar}).

Case (III)  $\lambda\in [\ell_k,\gamma_k)\cup(\gamma_k,\lambda_{k+1})$ for some $k\in \mathbb{N}$.
Take $\mathbf{i}=\lambda 0^k$ and $\mathbf{j}=0(1-\rho)^{k}$.
Note by (\ref{eq:composition-maps}) and (\ref{eq:lambda-gamma-eta}) that
\begin{equation*}\label{01-rhok=lambda0k}
  \begin{split}
 f_{0(1-\rho)^{k}}(0) =\rho(1-\rho^k)=\gamma_k \neq \lambda=f_{\lambda 0^k}(0).
  \end{split}
\end{equation*}
Then $f_{0(1-\rho)^{k}}\neq f_{\lambda 0^k}$.
 Furthermore,
 note by (\ref{eq:composition-maps}) that
 \begin{equation}\label{lambda0k=01-rhok}
  \begin{split}
  &f_{\lambda 0^k}(H)  =\left(\lambda +\rho^{k+1}(\rho+\lambda), \lambda +\rho^{k+1}(1-\rho-\lambda)\right),\\
 &f_{0(1-\rho)^{k}}(I)  =\left[\rho(1-\rho^{k}), \rho\right] .
  \end{split}
\end{equation}
 Since   $\lambda\in [\ell_k,\gamma_k)\cup(\gamma_k,\lambda_{k+1})$, by (\ref{eq:lambda-gamma-eta}) we have {$\lambda +\rho^{k+1}(\rho+\lambda)<\rho$},  and by Lemma \ref{inequ} (ii) it follows that
 $$
\rho(1-\rho^{k})=f_{0(1-\rho)^{k}}(0)<f_{\lambda 0^k}(1-\rho-\lambda)=\lambda +\rho^{k+1}(1-\rho-\lambda).
 $$
  Therefore, by using (\ref{lambda0k=01-rhok}) we obtain $f_{0(1-\rho)^{k}}(I)\cap f_{\lambda 0^k}(H)\neq \emptyset$, completing the proof.
\end{proof}

The proof for $\la\in\bigcup_{k=1}^\f\set{\la_k,\ga_k}$ to be sufficient is more involved. First we consider $\la=\la_k$ for some $k\in\N$.

\begin{proposition}\label{rho(1-rho)/(1+rho)}
  If $\lambda=\lambda_k$ for some $k\in \mathbb{N}$, then $E_\lambda$ is totally self-similar.
\end{proposition}
\begin{proof}Let $\la=\la_k =\frac{\rho(1-\rho^k)}{1+\rho^k}$.
 By Proposition \ref{equ state} it suffices to prove that for any $n\in \mathbb{N}$ and for any $\mathbf{i}, \mathbf{j} \in \Omega_{\la}^n$ with $f_{\mathbf{i}} \neq f_{\mathbf{j}}$ we have
 \begin{equation}\label{empty}
 f_{\mathbf{i}}(H) \cap f_{\mathbf{j}}(I)=\left(f_{\mathbf{i}}(0)+\rho^{n}(\rho+\la), f_{\mathbf{i}}(0)+\rho^{n}(1-\rho-\la)\right) \cap\left[f_{\mathbf{j}}(0), f_{\mathbf{j}}(0)+\rho^{n}\right]=\emptyset ,
 \end{equation}
which is equivalent to
\[f_{\mathbf{i}}(0)+\rho^{n}(\rho+\la)\geq f_{\mathbf{j}}(0)+\rho^{n} \quad \text{ or } \quad f_{\mathbf{i}}(0)+\rho^{n}(1-\rho-\la)\leq f_{\mathbf{j}}(0).\]
In other words,
\[\left|f_{\mathbf{i}}(0)-f_{\mathbf{j}}(0)\right|\geq \rho^n(1-\rho-\la).\]
So, we only need to prove that for any $n\in\mathbb{N}$ and for any $\mathbf{i}=i_1 \ldots i_n, \mathbf{j}=j_1 \ldots j_n \in \Omega_{\la}^n$ with $f_{\mathbf{i}} \neq f_{\mathbf{j}}$,
{\begin{align*}
  1-\rho-\la \leq  \frac{\left|f_{\mathbf{i}}(0)-f_{\mathbf{j}}(0)\right|}{\rho^{n}}&=\left|\sum_{m=1}^n \frac{\rho^{m-1}i_{m}}{\rho^{n}} -\sum_{m=1}^n \frac{\rho^{m-1}j_{m }}{\rho^{n}}\right|=\left|\sum_{m=1}^n \frac{i_m-j_m}{\rho^{n-m+1}} \right|=:\left|\sum_{m=1}^{n} \frac{c_m}{\rho^{m}} \right|,
\end{align*}}
where each {$c_m\in \Omega_{\la}^{\pm}=\left\{0, \pm \la, \pm\left(1-\rho-2\la\right),\pm\left(1-\rho-\la\right), \pm\left(1-\rho\right)\right\}$.}

Suppose on the contrary there exists a block $c_1c_2\cdots c_n\in (\Omega_{\la}^{\pm})^n$ for some $n\in \mathbb{N}$ such that
\begin{equation}\label{proof of contradiction}
  0<\left|\sum_{m=1}^{n} \frac{c_m}{\rho^{m}}\right|<1-\rho-\la.
 \end{equation}
 Furthermore, we can choose the block $c_1c_2\cdots c_n\in (\Omega_{\la}^{\pm})^n$ satisfying (\ref{proof of contradiction}) such that the finite sequence $\left|c_{n}\right|,\left|c_{n-1}\right|, \cdots,\left|c_1\right|$ is lexicographically minimal, and $c_n\ne 0$.
 Without loss of generality we can assume that $c_n>0$.
 Note that for any $\lambda\in(0,\rho)$,
 \begin{equation}\label{relationof4}
  1-\rho>1-\rho-\lambda>1-\rho-2\lambda>\lambda>0.
 \end{equation}
 If $n=1$, then by using $\la=\la_k\geq \lambda_1= \frac{\rho(1-\rho)}{1+\rho}$ and (\ref{relationof4}) it follows that
 \[\left|\frac{c_1}{\rho}\right|\ge \frac{\la}{\rho}\geq 1-\rho-\la,\]
 leading to a contradiction with
 (\ref{proof of contradiction}). So we must have $n\geq 2$.

 Note that $c_n\in\Om_\la^{\pm}$ and $c_n>0$. Then $c_n\in\set{\la, 1-\rho-2\la, 1-\rho-\la, 1-\rho}$.
If $c_n\in \left\{ 1-\rho-2\la, 1-\rho-\la, 1-\rho\right\}$, then by using $0<\la<\rho\le 1/4$ and $n\geq 2$ it follows that
$$\begin{aligned}
  \left|\sum_{m=1}^n \frac{c_m}{\rho^{m}}\right| &\geq\frac{1-\rho-2\la}{\rho^n} -\sum_{m=1}^{n-1}\frac{1-\rho}{\rho^m}=\frac{1-\rho-2\la}{\rho^n}-\frac{1-\rho^{n-1}}{\rho^{n-1}} \\
  &=1+\frac{1-2\rho-2\la}{\rho^n}>1,
  \end{aligned}$$
contradicting to ($\ref{proof of contradiction}$).
So, in the following it suffices to consider $c_n=\la$, which will be split into the following three  cases: {(I) $\la=\la_1$; (II) $\la=\la_k$ with $k\geq 2$ and $n\leq k$; (III) $\la=\la_k$ with $k\geq 2$ and $n>k$.}

Case (I) $c_n=\la=\lambda_1$.
If $c_{n-1}\leq 0$, then by using $c_{n-1}\in \Omega^{\pm}_{\lambda_1}$ and Lemma \ref{lambda1} we have  $c_{n-1}+\frac{\lambda_{1}}{\rho}\in \Omega_{\lambda_1}^{\pm}$.
Thus \[\sum_{m=1}^n \frac{c_m}{\rho^{m}}=\frac{1}{\rho^{n-1}}\left(c_{n-1}+\frac{\lambda_{1}}{\rho}\right)+\sum_{m=1}^{n-2} \frac{c_m}{\rho^{m}},\]
which contradicts to our assumption that $\left|c_n\right|,\left|c_{n-1}\right|, \cdots,\left|c_1\right|$ is lexicographically minimal.
So, $c_{n-1}>0$, and thus by using
  $0<\rho\le 1/4$ and (\ref{relationof4}) it follows that
$$\begin{aligned}
 \left|\sum_{m=1}^n \frac{c_m}{\rho^{m}}\right| &\geq\frac{\lambda}{\rho^n}+\frac{\lambda}{\rho^{n-1}} -\sum_{m=1}^{n-2}\frac{1-\rho}{\rho^m}=\frac{\lambda}{\rho^{n-1}}\left(\frac{1}{\rho}+1\right)-\frac{1-\rho^{n-2}}{\rho^{n-2}} \\
 &=1 +\frac{\la(1+\rho)-\rho^2}{\rho^{n}}>1,
 \end{aligned}$$
 where the last inequality follows by $\la=\la_1=\frac{\rho(1-\rho)}{1+\rho}$. This leads to a contradiction with
  (\ref{proof of contradiction}).

Case (II) $c_n=\la=\la_k$ with $k\geq 2$ and $n\leq k$.
Then by (\ref{eq:lambda-gamma-eta}) it follows that
 $$\begin{aligned}
  \left|\sum_{m=1}^n \frac{c_m}{\rho^{m}}\right| &\geq\frac{\la}{\rho^n}-\sum_{m=1}^{n-1}\frac{1-\rho}{\rho^m}=\frac{1-\rho^{k}}{\rho^{n-1}(1+\rho^{k})}-\frac{1-\rho^{n-1}}{\rho^{n-1}} \\
  &=1-\frac{2\rho^{k-n+1}}{1+\rho^k}\geq 1-\frac{2\rho}{1+\rho^k}=1-\rho-\la,
  \end{aligned}$$
  leading to a contradiction with (\ref{proof of contradiction}).

Case (III) $c_n= \la=\la_k$ with $k\geq 2$ and $n>k$.
We consider two subcases.

(III A)  $c_{n-1}c_{n-2}\cdots c_{n-k+1}$ has a digit $-(1-\rho-2\la)$.
Then by (\ref{eq:lambda-gamma-eta}) we have
\begin{equation*}
  \begin{aligned}
  \left|\sum_{m=1}^n \frac{c_m}{\rho^{m}}\right| &\geq\frac{\la}{\rho^n}-\sum_{m=n-k+2}^{n-1}\frac{1-\rho}{\rho^m}-\frac{1-\rho-2\la}{\rho^{n-k+1}}-\sum_{m=1}^{n-k}\frac{1-\rho}{\rho^m} \\
  &=\frac{\la}{\rho^n}+\frac{2\la}{\rho^{n-k+1}}-\sum_{m=1}^{n-1}\frac{1-\rho}{\rho^m}\\
  &=\frac{1-\rho^k}{\rho^{n-1}(1+\rho^k)}+\frac{2(1-\rho^k)}{\rho^{n-k}(1+\rho^k)}-\frac{1-\rho^{n-1}}{\rho^{n-1}}\\
&=1+\frac{2(1-\rho-\rho^k)}{\rho^{n-k}(1+\rho^k)}> 1,
  \end{aligned}
\end{equation*}
  contradicting to (\ref{proof of contradiction}).

(III B)  The digit  $-(1-\rho-2\la)$ does not appear in $c_{n-1}c_{n-2}\cdots c_{n-k+1}$.
Then we claim that at least one digit of $c_{n-1}c_{n-2} \cdots c_{n-k}$ should be positive.
Otherwise, by using $c_n=\la=\sum_{m=1}^{k-1} \rho^m(1-\rho)+\rho^k(1-\rho -\la)$ it follows that
\begin{equation}\label{c+1-rho-lambda}
\sum_{m=1}^n \frac{c_m}{\rho^{m}}=\sum_{m=n-k+1}^{n-1} \frac{c_m+1-\rho}{\rho^{m}}+\frac{c_{n-k}+1-\rho-\la}{\rho^{n-k}}+\sum_{m=1}^{n-k-1} \frac{c_m}{\rho^{m}}.
\end{equation}
Note that $c_m \in \Omega_{\la}^{\pm}$, $c_m\leq 0$ and $c_m\neq -(1-\rho-2\la)$ for all $n-k+1\leq m\leq n-1$.
Then $c_m+1-\rho\in \Omega_{\la}^{\pm}$.
Furthermore, since $c_{n-k}\in \Omega_{\la}^{\pm}$ and $c_{n-k}\leq 0$, we also have $c_{n-k}+1-\rho-\la\in \Omega_{\la}^{\pm}$.
Therefore, (\ref{c+1-rho-lambda}) gives another representation of $\sum_{m=1}^n \frac{c_m}{\rho^{m}}$, which is lexicographically smaller than $\left|c_n\right|,\left|c_{n-1}\right|, \cdots,\left|c_1\right|$, leading to a contradiction with our assumption.
By the claim, (\ref{relationof4}) and using $0<\rho\le 1/4$ it follows that $$\begin{aligned}
  \left|\sum_{m=1}^n \frac{c_m}{\rho^{m}}\right| &\geq\frac{\la}{\rho^n}-\sum_{m=n-k+1}^{n-1}\frac{1-\rho}{\rho^m}+\frac{\la}{\rho^{n-k}}-\sum_{m=1}^{n-k-1}\frac{1-\rho}{\rho^m}\\
  &=\frac{\la}{\rho^n}+\frac{\la+1-\rho}{\rho^{n-k}}-\sum_{m=1}^{n-1}\frac{1-\rho}{\rho^m}\\
  &=\frac{1-\rho^k}{\rho^{n-1}(1+\rho^k)}+\frac{\rho(1-\rho^k)+(1-\rho)(1+\rho^k)}{\rho^{n-k}(1+\rho^k)}-\frac{1-\rho^{n-1}}{\rho^{n-1}} \\
&=1+\frac{(1-2\rho)(1+\rho^k)}{\rho^{n-k}(1+\rho^k)}> 1,
 \end{aligned}$$
 contradicting to (\ref{proof of contradiction}).

 Hence, by Cases (I)--(III)  {it} follows that (\ref{proof of contradiction}) fails, and then  proves (\ref{empty}) as required.
\end{proof}

\begin{proof}[Proof of Theorem \ref{th:totally-self-similar} (i)]
  Let $0<\lambda<\rho\le 1/4$. By Propositions \ref{NTSS} and \ref{rho(1-rho)/(1+rho)}  we only need to prove that if $\lambda=\gamma_{k}=\rho(1-\rho^k)$ for some $k\in \mathbb{N}$, then $E_\lambda$ is totally self-similar.
Take $\la=\ga_k$ with $k\in\N$. By the same argument as in the proof of Proposition \ref{rho(1-rho)/(1+rho)} it suffices to prove that for any $n\in \mathbb{N}$ and for any {$c_m\in \Omega_{\la}^{\pm}=\left\{0, \pm \la, \pm\left(1-\rho-2\la\right),\pm\left(1-\rho-\la\right),\pm\left(1-\rho\right)\right\} $} with $1\leq m\leq n$ we have
$|\sum_{m=1}^\f\frac{c_m}{\rho^m}|\ge 1-\rho-\la$. Indeed, we can prove
 \[\left|\sum_{m=1}^{n} \frac{c_m}{\rho^{m}} \right|\geq 1-\rho. \]

   Suppose on the contrary there exists a block $c_1c_2\cdots c_n\in (\Omega_{\la}^{\pm})^n$ with $n\in \mathbb{N}$ such that
 \begin{equation}\label{proof of contradiction2}
   0<\left|\sum_{m=1}^{n} \frac{c_m}{\rho^{m}}\right|<1-\rho.
  \end{equation}
Furthermore, we can choose the block $c_1c_2\cdots c_n\in (\Omega_{\la}^{\pm})^n$ satisfying (\ref{proof of contradiction2}) such that the finite sequence $\left|c_{n}\right|,\left|c_{n-1}\right|, \cdots,\left|c_1\right|$ is lexicographically minimal and $c_n\neq 0$.
Without loss of generality we can assume that $c_{n}>0$.
Note that $\la=\ga_k\geq  \rho(1-\rho)$. Then by (\ref{proof of contradiction2}) and the same argument as in the proof of Proposition \ref{rho(1-rho)/(1+rho)} we have $n\geq 2$.

Note that $c_n\in\Om_\la^{\pm}$ and $c_n>0$. Then $c_n\in\set{\la, 1-\rho-2\la, 1-\rho-\la, 1-\rho}$. If $c_n\in \left\{1-\rho-2\la, 1-\rho-\la, 1-\rho\right\}$, then by using $0<\la<\rho\le1/4$ and (\ref{relationof4}) it follows that
$$\begin{aligned}
\left|\sum_{m=1}^{n} \frac{c_m}{\rho^{m}}\right| &\geq\frac{1-\rho-2\la}{\rho^{n}} -\sum_{m=1}^{n-1}\frac{1-\rho}{\rho^m}
 =1+\frac{1-2\rho-2\la}{\rho^n} >1,
\end{aligned}$$
leading to a contradiction with (\ref{proof of contradiction2}).
%
 {So, in the following it suffices to consider $c_n=\la$, which will be split into the following two cases: (I) $n\leq k$; (II)  $n>k$.}

{Case (I) $c_n=\la$ with $n\leq k$.
Then by using $\la=\ga_k=\rho(1-\rho^k)$  it follows that
 $$\begin{aligned}
  \left|\sum_{m=1}^n \frac{c_m}{\rho^{m}}\right| &\geq\frac{\la}{\rho^n}-\sum_{m=1}^{n-1}\frac{1-\rho}{\rho^m}=\frac{1-\rho^{k}}{\rho^{n-1}}-\frac{1-\rho^{n-1}}{\rho^{n-1}} \\
  &=1-\rho^{k-n+1}\geq 1-\rho,
  \end{aligned}$$
  leading to a contradiction with (\ref{proof of contradiction2}).}

 {Case (II) $c_n= \la$ with $n>k$.
We consider two subcases.}

 {(II A) $c_{n-1}c_{n-2}\cdots c_{n-k}$ contains a digit $-(1-\rho-2\la)$.
Then by (\ref{relationof4}) and $\la=\ga_k=\rho(1-\rho^k)$ it follows that
\begin{equation}\label{gamma-to+}
  \begin{aligned}
  \left|\sum_{m=1}^{n} \frac{c_m}{\rho^{m}}\right| &\geq\frac{\la}{\rho^{n}}-\sum_{m=n-k+1}^{n-1}\frac{1-\rho}{\rho^m}-\frac{1-\rho-2\la}{\rho^{n-k}}-\sum_{m=1}^{n-k-1}\frac{1-\rho}{\rho^m} \\
  &=\frac{\la}{\rho^{n}}+\frac{2\la}{\rho^{n-k}}-\sum_{m=1}^{n-1}\frac{1-\rho}{\rho^m}\\
   &=1+\frac{1-2\rho^{k}}{\rho^{n-k-1}} \ge 2-2\rho^{k}>1-\rho,
  \end{aligned}
\end{equation}
 where the last inequality follows by $0<\rho\le1/4$. This leads to a contradiction with (\ref{proof of contradiction2}).}

(II B)  The digit  $-(1-\rho-2\la)$ never occurs in the block $c_{n-1}c_{n-2}\cdots c_{n-k}$.
  Then we claim that at least one digit in $c_{n-1}c_{n-2} \cdots c_{n-k}$ is positive.
  Otherwise, by using $c_n=\la=\rho-\rho^{k+1}=\sum_{m=1}^{k} \rho^m(1-\rho)$ we have
  \begin{equation}\label{c+1-rho-lambda3}
  \sum_{m=1}^{n} \frac{c_m}{\rho^{m}}=\sum_{m=n-k}^{n-1} \frac{c_m+1-\rho}{\rho^{m}}+\sum_{m=1}^{n-k-1} \frac{c_m}{\rho^{m}}.
  \end{equation}
  Note that $c_{m}\in\Omega_{\la}^{\pm}$, $c_{m}\leq 0$ and $c_{m}\neq -(1-\rho-2\la)$ for any $n-k\leq m\leq n-1$.
  Then $c_{m}+1-\rho \in \Omega_{\la}^{\pm}$ for all $n-k\leq m\leq n-1$.
  Therefore, (\ref{c+1-rho-lambda3}) contradicts to the minimality of $\left|c_{n}\right|,\left|c_{n-1}\right|, \cdots,\left|c_1\right|$. This proves the claim.
So by (\ref{relationof4}) it follows that
  $$\begin{aligned}
  \left|\sum_{m=1}^n \frac{c_m}{\rho^{m}}\right| &\geq\frac{\la}{\rho^n}-\sum_{m=n-k+1}^{n-1}\frac{1-\rho}{\rho^m}+\frac{\la}{\rho^{n-k}}-\sum_{m=1}^{n-k-1}\frac{1-\rho}{\rho^m}\\
  &>\frac{\la}{\rho^{n}}-\sum_{m=n-k+1}^{n-1}\frac{1-\rho}{\rho^m}-\frac{1-\rho-2\la}{\rho^{n-k}}-\sum_{m=1}^{n-k-1}\frac{1-\rho}{\rho^m} > 1-\rho,
 \end{aligned}$$
 where the last inequality holds by the same argument as in   (\ref{gamma-to+}). Again this leads to a contradiction with
  (\ref{proof of contradiction2}).

  Therefore, (\ref{proof of contradiction2}) fails by Cases (I) and (II). This completes the proof.
\end{proof}


\subsection{{Total self-similarity}  of $E_\lambda$ for $\lambda\in \left(\frac{1-2\rho}{2},\frac{1-\rho}{2}\right)$}
The proof of Theorem \ref{th:totally-self-similar} (ii) is similar to that for Theorem \ref{th:totally-self-similar} (i).
Take $\lambda\in \left(\frac{1-2\rho}{2},\frac{1-\rho}{2}\right)$. Then the   hole $H$ is given by (see the right graph of Figure \ref{Fig:1})
  \[H=I \setminus \bigcup_{d\in\Om_\la}f_{d}(I)=(\rho,\lambda)\cup(1-\lambda,1-\rho).\]
 Let $H_1:=(\rho,\lambda)$ and $H_2:=(1-\lambda,1-\rho)$. Then $H=H_1\cup H_2$ with the {union disjoint}.
  Recall from (\ref{eq:lambda-gamma-eta}) that $\eta_k=\frac{1-2\rho+\rho^{k+1}}{2}\in \left(\frac{1-2\rho}{2},\frac{1-\rho}{2}\right)$ for all $k\in\N$, and $\eta_k\searrow \frac{1-2\rho}{2}$ as $k\to \infty$.
\begin{proposition}\label{(1-2rho+rho(k+1))/2}
  If $\lambda=\eta_k$ for some $k\in \mathbb{N}$, then $E_\lambda$ is totally self-similar.
\end{proposition}

\begin{proof}
Take $\la=\eta_k$.
By Proposition \ref{equ state} it suffices to prove that for any $n\in \mathbb{N}$ and for any $\mathbf{i}, \mathbf{j} \in \Omega_{\lambda}^n$ with $f_{\mathbf{i}} \neq f_{\mathbf{j}}$ we have
\begin{align*}\label{empty2}
f_{\mathbf{i}}(H) \cap f_{\mathbf{j}}(I)&=f_{\mathbf{i}}(H_1\cup H_2) \cap f_{\mathbf{j}}(I)=\left(f_{\mathbf{i}}(H_1) \cap f_{\mathbf{j}}(I)\right)\cup \left(f_{\mathbf{i}}(H_2) \cap f_{\mathbf{j}}(I)\right)\\
&=\left(\left(f_{\mathbf{i}}(0)+\rho^{n+1}, f_{\mathbf{i}}(0)+\rho^{n}\la \right) \cap\left[f_{\mathbf{j}}(0), f_{\mathbf{j}}(0)+\rho^{n}\right]\right)\\
&\quad \cup\left(\left(f_{\mathbf{i}}(0)+\rho^{n}(1-\la), f_{\mathbf{i}}(0)+\rho^{n}(1-\rho)\right) \cap\left[f_{\mathbf{j}}(0), f_{\mathbf{j}}(0)+\rho^{n}\right]\right)\\
&=\emptyset ,
\end{align*}
which is equivalent to
\[\left|f_{\mathbf{i}}(0)-f_{\mathbf{j}}(0)\right|\geq \max\left\{\rho^n \la, \rho^n(1-\rho)\right\}=\rho^n(1-\rho),\]
where the equality holds since $\la<\frac{1-\rho}{2}$.
So, we only need to prove that for any $n\in\mathbb{N}$ and for any $\mathbf{i}=i_1 \ldots i_n, \mathbf{j}=j_1 \ldots j_n \in \Omega_{\la}^n$ with $f_{\mathbf{i}} \neq f_{\mathbf{j}}$,
{\begin{align*}
  1-\rho\leq  \frac{\left|f_{\mathbf{i}}(0)-f_{\mathbf{j}}(0)\right|}{\rho^{n}}&=\left|\sum_{m=1}^n \frac{\rho^{m-1}i_{m}}{\rho^{n}} -\sum_{m=1}^n \frac{\rho^{m-1}j_{m }}{\rho^{n}}\right|=\left|\sum_{m=1}^n \frac{i_m-j_m}{\rho^{n-m+1}} \right|=:\left|\sum_{m=1}^{n} \frac{c_m}{\rho^{m}} \right|,
 \end{align*}}
where each {$c_m\in \Omega_{\la}^{\pm}=\left\{0, \pm \la,\pm\left(1-\rho-2\la\right),\pm\left(1-\rho-\la\right), \pm\left(1-\rho\right)\right\}$}.

Suppose on the contrary there exists a block $c_1c_2\cdots c_n\in (\Omega_{\la}^{\pm})^n$ for some $n\in \mathbb{N}$ such that
\begin{equation}\label{proof of contradiction3}
0<\left|\sum_{m=1}^{n} \frac{c_m}{\rho^{m}}\right|<1-\rho.
\end{equation}
By the same argument as in the proof of Proposition \ref{rho(1-rho)/(1+rho)}, we can choose the block $c_1c_2\cdots c_n\in (\Omega_{\la}^{\pm})^n$ satisfying (\ref{proof of contradiction3}) such that the finite sequence $\left|c_{n}\right|,\left|c_{n-1}\right|, \cdots,\left|c_1\right|$ is lexicographically minimal and $c_n> 0$.
Note that  $\la=\eta_k\le \eta_1=\frac{1-2\rho+\rho^2}{2}$ and $\rho\in(0,1/4]$. Then  by (\ref{proof of contradiction3}) it follows that $n\geq 2$.
Furthermore, by using $0<\rho\le1/4$ and $\frac{1-2\rho}{2}<\la<\frac{1-\rho}{2}$ we have
\begin{equation}\label{eta5}
  0<1-\rho-2\la<\rho\le\la<1-\rho-\la<1-\rho .
\end{equation}
Note that $c_n\in\Om_\la^{\pm}$ and $c_n>0$. Then $c_n\in\set{1-\rho-2\la, \la, 1-\rho-\la, 1-\rho}$.
If $c_n\in \left\{\la, 1-\rho-\la, 1-\rho\right\}$, then by (\ref{eta5}) we have
$$\begin{aligned}
  \left|\sum_{m=1}^{n} \frac{c_m}{\rho^{m}}\right| &\geq\frac{\la}{\rho^{n}} -\sum_{m=1}^{n-1}\frac{1-\rho}{\rho^m}=\frac{\la}{\rho^{n}}-\frac{1-\rho^{n-1}}{\rho^{n-1}}=1+\frac{\la-\rho}{\rho^n}\geq 1,
  \end{aligned}$$
leading to a contradiction with (\ref{proof of contradiction3}).
So, in the following it suffices to consider $c_n=1-\rho-2\la$, which will be split into the following two cases: {(I)   $n\leq k$; (II) $n>k$}.

  {Case (I) $c_n=1-\rho-2\la$ with $n\leq k$.
Then by (\ref{eta5}) and using $\lambda=\eta_k=\frac{1-2\rho+\rho^{k+1}}{2}$ it follows that
 $$\begin{aligned}
  \left|\sum_{m=1}^n \frac{c_m}{\rho^{m}}\right| &\geq\frac{1-\rho-2\la}{\rho^n}-\sum_{m=1}^{n-1}\frac{1-\rho}{\rho^m}=\frac{1-\rho^{k}}{\rho^{n-1}}-\frac{1-\rho^{n-1}}{\rho^{n-1}} \\
  &=1-\rho^{k-n+1}\geq 1-\rho,
  \end{aligned}$$
  leading to a contradiction with (\ref{proof of contradiction3}).}

 {Case (II) $c_n= 1-\rho-2\la$ with $n>k$.
We consider two subcases.}

 (II A) $c_{n-1}c_{n-2}\cdots c_{n-k}$ contains a digit $-(1-\rho-2\la)$.
Then by (\ref{eta5}) and $\lambda=\eta_k=\frac{1-2\rho+\rho^{k+1}}{2}$ it follows that
\begin{equation}\label{eta-to+}
  \begin{aligned}
  \left|\sum_{m=1}^{n} \frac{c_m}{\rho^{m}}\right| &\geq\frac{1-\rho-2\la}{\rho^{n}}-\sum_{m=n-k+1}^{n-1}\frac{1-\rho}{\rho^m}-\frac{1-\rho-2\la}{\rho^{n-k}}-\sum_{m=1}^{n-k-1}\frac{1-\rho}{\rho^m} \\
  &=\frac{1-\rho-2\la}{\rho^{n}}+\frac{2\la}{\rho^{n-k}}-\sum_{m=1}^{n-1}\frac{1-\rho}{\rho^m}\\
  &=1+\frac{1-3\rho+\rho^{k+1}}{\rho^{n-k}} >1,
  \end{aligned}
\end{equation}
 where the last inequality follows by $0<\rho\le1/4$. This leads to a contradiction with (\ref{proof of contradiction3}).

 {(II B)} The digit  $-(1-\rho-2\la)$ never occurs in $c_{n-1}c_{n-2}\cdots c_{n-k}$.
Then we claim that at least one digit of $c_{n-1}c_{n-2} \cdots c_{n-k}$ is positive.
Otherwise, by using $c_n=1-\rho-2\la=\rho-\rho^{k+1}=\sum_{m=1}^{k} \rho^m(1-\rho)$ we have
\begin{equation}\label{c+1-rho-lambda2}
\sum_{m=1}^{n} \frac{c_m}{\rho^{m}}=\sum_{m=n-k}^{n-1} \frac{c_m+1-\rho}{\rho^{m}}+\sum_{m=1}^{n-k-1} \frac{c_m}{\rho^{m}}.
\end{equation}
Note that $c_{m}\in\Omega_{\la}^{\pm}$, $c_{m}\leq 0$ and $c_{m}\neq -(1-\rho-2\la)$ for any $n-k\leq m\leq n-1$.
Then $c_{m}+1-\rho \in \Omega_{\la}^{\pm}$ for all $n-k\le m\le n-1$.
Therefore, (\ref{c+1-rho-lambda2}) contradicts to the minimality of $\left|c_{n}\right|,\left|c_{n-1}\right|, \cdots,\left|c_1\right|$.
This proves the claim, and then by (\ref{eta5})
 it follows that
 $$\begin{aligned}
 \left|\sum_{m=1}^n \frac{c_m}{\rho^{m}}\right| &\geq\frac{1-\rho-2\la}{\rho^n}-\sum_{m=n-k+1}^{n-1}\frac{1-\rho}{\rho^m}+\frac{1-\rho-2\la}{\rho^{n-k}}-\sum_{m=1}^{n-k-1}\frac{1-\rho}{\rho^m}\\
 &\geq\frac{1-\rho-2\la}{\rho^{n}}-\sum_{m=n-k+1}^{n-1}\frac{1-\rho}{\rho^m}-\frac{1-\rho-2\la}{\rho^{n-k}}-\sum_{m=1}^{n-k-1}\frac{1-\rho}{\rho^m} > 1,
\end{aligned}$$
where the last inequality holds by the same argument as in   (\ref{eta-to+}).
This again contradicts  to (\ref{proof of contradiction3}).

 By Cases (I) and (II) it follows that (\ref{proof of contradiction3}) does not hold, and thus $f_{\mathbf{i}}(H) \cap f_{\mathbf{j}}(I)= \emptyset $ for any $\i, \j\in\Om_\la^n$ with $n\in\N$. This  completes the proof.
\end{proof}


\begin{proof}[Proof of Theorem \ref{th:totally-self-similar} (ii)]

By Proposition \ref{(1-2rho+rho(k+1))/2}  we only need to prove that if $\lambda\in (\frac{1-2\rho}{2},\frac{1-\rho}{2})\setminus \bigcup_{k=1}^{\infty} \left\{ \eta_k\right\}$, then $E_\lambda$ is not totally self-similar.
Note that $\eta_0=\frac{1-\rho}{2}$ and $\eta _{k}\searrow \frac{1-2\rho}{2}$ as $k\to \infty$.
So by Proposition \ref{equ state} it suffices to show that for any $k\in \mathbb{N}$ and for any $\lambda\in (\eta_k,\eta_{k-1})$ we can find two words $\i, \j\in\Om_\la^n$ with $n\in\N$ such that
\begin{equation}\label{eq:Not-totally-self-similar-1}
f_\i\ne f_\j\quad\textrm{and}\quad f_\i(I)\cap f_\j(H)\ne\emptyset,
\end{equation} where $I=[0,1]$ and $H=(\rho,\la)\cup(1-\la, 1-\rho)$.

Let $\la\in(\eta_k,\eta_{k-1})$  {for some $k\in \mathbb{N}$}, and take $\mathbf{i}=\lambda(1-\rho)^{k-1}, \mathbf{j}=(1-\rho-\lambda) 0^{k-1}\in\Om_\la^k$.
Since $\lambda<\eta_{k-1}=\frac{1-2\rho+\rho^k}{2}$, by Lemma (\ref{inequ}) (iii) we have
\begin{align*}
  f_{\lambda(1-\rho)^{k-1}}(0)<f_{(1-\rho-\lambda) 0^{k-1}}(0).
 \end{align*}
Then $f_{\lambda(1-\rho)^{k-1}}\neq f_{(1-\rho-\lambda) 0^{k-1}}$.
Furthermore, since $H_1=(\rho,\lambda)\subset H$, by (\ref{eq:composition-maps}) it follows that
 \begin{equation}\label{l-ambda-rho0k=lambda1-rhok}
 \begin{split}
f_{\lambda(1-\rho)^{k-1}}(I) &=\left[\rho-\rho^{k}+\lambda, \rho+\lambda\right],\\
 f_{(1-\rho-\lambda) 0^{k-1}}(H_1) &=\left(1-\rho-\lambda+\rho^{k+1},1-\rho-\lambda+\rho^{k} \lambda \right).
 \end{split}
\end{equation}
Note that $\lambda\in  (\eta_k,\eta_{k-1})$. Then $1-\rho-\lambda+\rho^{k+1}<\rho+\lambda$.
Furthermore, by Lemma \ref{inequ} (iii) we obtain
$$
\rho-\rho^{k}+\lambda=f_{\lambda(1-\rho)^{k-1}}(0)<f_{(1-\rho-\lambda) 0^{k-1}}(0)<f_{(1-\rho-\lambda) 0^{k-1}}(\lambda)=1-\rho-\lambda+\rho^{k}\lambda.
$$
So, by (\ref{l-ambda-rho0k=lambda1-rhok}) it follows that  {$f_{\lambda(1-\rho)^{k-1}}(I)\cap f_{(1-\rho-\lambda) 0^{k-1}}(H_1)\neq \emptyset $}, proving (\ref{eq:Not-totally-self-similar-1}).
\end{proof}

\section{The spectrum of $E_\lambda$}\label{sec:spectrum}
Recall by Definition \ref{def of spec} that the spectrum $l_\la$ of $E_\lambda$ is given by
$$
 l_\lambda=\inf\left\{\frac{\left|f_{\mathbf{i}}(0)-f_{\mathbf{j}}(0)\right|}{\rho ^{n}}: \mathbf{i}, \mathbf{j} \in \Omega_\lambda^n \text { with } f_{\mathbf{i}} \neq f_{\mathbf{j}};~  n\in \mathbb{N}\right\}.
$$In this section we will characterize the total self-similarity of $E_\la$ by using the spectrum, and prove Theorem \ref{spectrum}.
First we consider the spectrum $l_\lambda$ when $E_\lambda$ is totally self-similar.

\begin{proposition}\label{spectrum of tss}\mbox{}
  \begin{enumerate}[\rm(i)]
 \item If $\lambda=\lambda_k$ for some $k\in\mathbb{N}$, then $l_{\lambda} = 1-\rho-\lambda$.
 \item If $\lambda=\gamma_k$   for some $k\in\mathbb{N}$, then $l_{\lambda}= 1-\rho$.
  \item If $\lambda= \eta_k$ for some $k\in\mathbb{N}$, then $l_{\lambda}= 1-\rho$.
  \end{enumerate}
\end{proposition}
\begin{proof}

For (i) let $\lambda=\lambda_k=\frac{ \rho(1-\rho ^k)}{1+\rho ^k}$ for some $k\in \mathbb{N}$.
  Then by the proof of Proposition \ref{rho(1-rho)/(1+rho)} it follows that
$$
l_{\la}=\inf \left\{\frac{\left|f_{\mathbf{i}}(0)-f_{\mathbf{j}}(0)\right|}{\rho ^n}: \mathbf{i}, \mathbf{j} \in \Omega_{\la}^n \text { with } f_{\mathbf{i}} \neq f_{\mathbf{j}} , n\in \mathbb{N}\right\} \geq 1-\rho-\la .
$$
On the other hand, let $\mathbf{i}=0(1-\rho)^{k-1}$ and $\mathbf{j}=\la 0^{k-1}$.
By using  {$\rho\in (0,1/4]$} and $\la=\frac{ \rho(1-\rho ^k)}{1+\rho ^k}$ it follows  by (\ref{eq:composition-maps}) that
\begin{equation*}
 f_{0(1-\rho)^{k-1}}(0)=\rho(1-\rho^{k-1})< \la=f_{\la 0^{k-1}}(0),
\end{equation*}
which implies $f_{0(1-\rho)^{k-1}}\neq f_{\la 0^{k-1}}$.
Then
 \begin{align*}
  \frac{\left|f_{0(1-\rho)^{k-1}}(0)-f_{\la 0^{k-1}}(0)\right|}{\rho ^k}&=\frac{\left|\rho(1-\rho^{k-1})-\la\right|}{\rho ^k} =\frac{1+\rho ^k-2\rho}{1+\rho ^k}=1-\rho-\la.
  \end{align*}
This proves $l_{\la}=1-\rho-\la$.

Next we consider (ii). Let $\lambda=\gamma_k= \rho(1-\rho ^k)$ for some $k\in \mathbb{N}$.
Then by the proof of Theorem \ref{th:totally-self-similar} (i) it follows that
$$
l_{\la}=\inf \left\{\frac{\left|f_{\mathbf{i}}(0)-f_{\mathbf{j}}(0)\right|}{\rho ^n}: \mathbf{i}, \mathbf{j} \in \Omega_{\la}^n \text { with } f_{\mathbf{i}} \neq f_{\mathbf{j}} , n\in \mathbb{N}\right\} \geq 1-\rho.
$$
On the other hand, take $\mathbf{i}=\la 0^{k-1}$ and $\mathbf{j}=0(1-\rho)^{k-1}$.
Then by using $\la=\rho(1-\rho^k)$ we have
\begin{equation*}
f_{\la 0^{k-1}}(0)=\la >\rho(1-\rho^{k-1})=f_{0(1-\rho)^{k-1}}(0),
\end{equation*}
which yields $f_{\la 0^{k-1}}\neq f_{0(1-\rho)^{k-1}}$. Furthermore,
\begin{align*}
\frac{\left|f_{\la 0^{k-1}}(0)-f_{0(1-\rho)^{k-1}}(0)\right|}{\rho ^k}&=\frac{\left|\la -\rho+\rho^{k}\right|}{\rho ^k}=\frac{\left|\rho-\rho^{k+1}-\rho+\rho^{k}\right|}{\rho ^k}=1-\rho.
\end{align*}
Thus, $l_{\la}=1-\rho$.

Finally we prove (iii). Let $\lambda=\eta_k= \frac{1-2\rho+\rho^{k+1}}{2}$ for some $k\in \mathbb{N}$.
Then by the proof of  {Proposition \ref{(1-2rho+rho(k+1))/2}} it follows that
$$
l_{\la}=\inf \left\{\frac{\left|f_{\mathbf{i}}(0)-f_{\mathbf{j}}(0)\right|}{\rho ^n}: \mathbf{i}, \mathbf{j} \in \Omega_{\la}^n \text { with } f_{\mathbf{i}} \neq f_{\mathbf{j}} , n\in \mathbb{N}\right\} \geq 1-\rho.
$$
On the other hand, let $\mathbf{i}=(1-\rho-\la) 0^{k-1}$ and $\mathbf{j}=\la(1-\rho)^{k-1}$.
By using $\la=\frac{1-2\rho+\rho^{k+1}}{2}$ it follows that
\begin{equation*}
 f_{(1-\rho-\la) 0^{k-1}}(0)=1-\rho-\la >\la+\rho-\rho^{k} =f_{\la(1-\rho)^{k-1}}(0),
\end{equation*}which implies $f_{(1-\rho-\la) 0^{k-1}}\ne f_{\la(1-\rho)^{k-1}}$.
Furthermore,
\begin{align*}
\frac{\left|f_{(1-\rho-\la) 0^{k-1}}(0)-f_{\la(1-\rho)^{k-1}}(0)\right|}{\rho ^k}&=\frac{\left|\la+\rho-\rho^{k}-1+\rho+\la\right|}{\rho ^k}\\
&=\frac{\left|2\rho-\rho^{k}-1+1-2\rho+\rho^{k+1}\right|}{\rho ^k}=1-\rho.
\end{align*}
So, $l_{\la}=1-\rho$.
\end{proof}

\begin{proof}[Proof of Theorem \ref{spectrum}]
By Proposition \ref{spectrum of tss} it suffices to prove (ii) and (iii).
First we prove (ii). By Theorem \ref{th:totally-self-similar} and  Proposition \ref{spectrum of tss}  it suffices to prove that for any $\la\in(0,\rho)\setminus\bigcup_{k=1}^\f\set{\la_k, \ga_k}$ we have $l_\la<1-\rho-\la$.
  Note by (\ref{eq:lambda-gamma-eta}) and (\ref{ell}) that $\lambda_k=\frac{\rho(1-\rho^{k})}{1+\rho^{k}}$, $\ell_k=\frac{\rho(1-\rho^{k+1})}{1+\rho^k+\rho^{k+1}}$, $\gamma_k=\rho(1-\rho^{k})$.
  Then $0<\lambda_{k}<\ell_{k}<\gamma_{k}<\lambda_{k+1}$ for any $k\in \mathbb{N}$, and $\la_k\nearrow \rho$ as $k\to\f$.
  So, it suffices to prove  $l_\lambda<1-\rho-\lambda$ in the following four cases:
  {(I)} $\lambda\in (0,\lambda_1)$;
  {(II)} $\lambda\in \bigcup_{k=1}^{\infty} (\lambda_k,\ell_k)$;
  {(III)} $\lambda\in \bigcup_{k=1}^{\infty} [\ell_k,\gamma_k)$;
  {(IV)} $\lambda\in \bigcup_{k=1}^{\infty} (\gamma_k,\lambda_{k+1})$.

  Case {(I)} $\lambda\in (0,\lambda_1)$.
  Take $\mathbf{i}=0$ and $\mathbf{j}=\lambda$.
  Since $\lambda<\lambda_1=\frac{\rho(1-\rho)}{1+\rho}$, we have $\frac{\lambda}{\rho}<1-\rho-\lambda$, and then
  \[  l_\lambda\leq \frac{\left|f_{0}(0)-f_{\lambda}(0)\right|}{\rho}=\frac{\lambda}{\rho}<1-\rho-\lambda.\]

  Case {(II)} $\lambda\in  (\lambda_k,\ell_k)$ for some $k\in \mathbb{N}$.
Take $\mathbf{i}=\lambda 0^k$ and $\mathbf{j}=0(1-\rho)^{k-1}(1-\rho-\lambda )$.
   Since $\lambda>\lambda_k=\frac{\rho(1-\rho^k)}{1+\rho^k}$, by (\ref{eq:composition-maps}) we have
  \begin{equation}\label{lambda0k01rhok-11-rho-lambda}
  f_{\lambda 0^k}(0)=\lambda>\rho-\rho^{k}\lambda-\rho^{k+1}=f_{0(1-\rho)^{k-1}(1-\rho-\lambda )}(0),
  \end{equation}
  which implies $ f_{\lambda 0^k}\ne f_{0(1-\rho)^{k-1}(1-\rho-\lambda )}$.
  Note that $\lambda<\ell_k=\frac{\rho(1-\rho^{k+1})}{1+\rho^k+\rho^{k+1}}$. Then   by (\ref{lambda0k01rhok-11-rho-lambda}) it follows that
\[  l_\lambda\leq \frac{\left|f_{\lambda 0^k}(0)-f_{0(1-\rho)^{k-1}(1-\rho-\lambda )}(0)\right|}{\rho^{k+1}}=\frac{\lambda+\rho^{k}\lambda-\rho+\rho^{k+1}}{\rho^{k+1}}<1-\rho-\lambda.\]

  Case {(III)}  $\lambda\in [\ell_k,\gamma_k)$ for some $k\in \mathbb{N}$.
  Let $\mathbf{i}=\lambda 0^k$ and $\mathbf{j}=0(1-\rho)^{k}$.
 Then
  \begin{align*}
 f_{0(1-\rho)^{k}}(0)=\rho(1-\rho^k)=\gamma_k>\lambda=f_{\lambda 0^k}(0),
  \end{align*}which implies
  $f_{0(1-\rho)^{k}}\neq f_{\lambda 0^k}$.
Furthermore,
\[  l_\lambda\leq \frac{\left|f_{\lambda 0^k}(0)-f_{0(1-\rho)^{k}}(0)\right|}{\rho^{k+1}}=\frac{\rho(1-\rho^k)-\lambda}{\rho^{k+1}}<1-\rho-\lambda,\]
where the last inequality follows by $0<\rho\le1/4$   that
\[\frac{\rho-2\rho^{k+1}+\rho^{k+2}}{1-\rho^{k+1}}<\frac{\rho(1-\rho^{k+1})}{1+\rho^k+\rho^{k+1}}=\ell_k\leq \lambda.\]

Case {(IV)} $\lambda\in (\gamma_k,\lambda_{k+1})$ for some $k\in \mathbb{N}$.
 Take $\mathbf{i}=\lambda 0^k$ and $\mathbf{j}=0(1-\rho)^{k}$.
 Similar to Case {(III)} we have $f_{0(1-\rho)^{k}}\neq f_{\lambda 0^k}$.
Since $\ga_k<\la<\la_{k+1}$, by (\ref{eq:lambda-gamma-eta}) it follows that
\[ 0<l_\lambda\leq \frac{\left|f_{\lambda 0^k}(0)-f_{0(1-\rho)^{k}}(0)\right|}{\rho^{k+1}}=\frac{\lambda-\rho(1-\rho^k)}{\rho^{k+1}}<1-\rho-\lambda.\]
By Cases {(I)--(IV)} we prove (ii).

Next we prove (iii). By Theorem \ref{th:totally-self-similar} and Proposition \ref{spectrum of tss} it suffices to prove that for any $\lambda \in \left(\frac{1-2\rho}{2},\frac{1-\rho}{2}\right)\setminus\bigcup_{k=1}^\f\set{\eta_k}$ we have $l_\la<1-\rho$. Note that $\eta_0=\frac{1-\rho}{2}$ and $\eta_k\searrow \frac{1-2\rho}{2}$ as $k\to\f$. Then we only need to prove $l_\la<1-\rho$ for any $\la\in(\eta_k, \eta_{k-1})$ with $k\in\N$. Now take
 $\lambda\in (\eta_k,\eta_{k-1})$ for some $k\in\N$.
Let $\mathbf{i}=(1-\rho-\lambda) 0^{k-1}$ and $\mathbf{j}=\lambda(1-\rho)^{k-1}$.
Since $\lambda<\eta_{k-1}$, by Lemma  \ref{inequ}  (iii) we have $f_{\lambda(1-\rho)^{k-1}}(0)<f_{(1-\rho-\lambda) 0^{k-1}}(0)$, which gives $f_{\lambda(1-\rho)^{k-1}}\neq f_{(1-\rho-\lambda) 0^{k-1}}$.
Furthermore, by using  {$\lambda>\eta_k=\frac{1-2\rho+\rho^{k+1}}{2}$}   it follows that
\[  l_\lambda\leq \frac{\left|f_{\lambda(1-\rho)^{k-1}}(0)-f_{(1-\rho-\lambda) 0^{k-1}}(0)\right|}{\rho^{k}}=\frac{1-\rho-\lambda-\rho+\rho^{k}-\lambda}{\rho^{k}}<1-\rho \]
as desired. This completes the proof.
\end{proof}

\section{Unique codings of $E_\la$ when $E_\la$ is totally self-similar}\label{sec:unique-codings-totally-self-similar}
Recall that $U_\la$ consists of all $x\in E_\la$ having a unique coding with respect to the IFS $\mathcal F_\la$ defined in (\ref{eq:IFS}).
In this section we will study the set $U_\la$  when $E_\la$ is totally self-similar, and prove Theorem \ref{th:unique-coding-totally-self-similar}. Let $\la\in(0,\rho)\cup(\frac{1-2\rho}{2}, \frac{1-\rho}{2})$. Note by Theorem \ref{th:totally-self-similar} that $E_\la$ is totally self-similar if and only if $\la\in\bigcup_{k=1}^\f\set{\la_k, \ga_k,\eta_k}$. So we will calculate the Hausdorff dimension of $U_\la$ for $\la=\la_k, \ga_k$ and $\eta_k$.

First we prove Theorem \ref{th:unique-coding-totally-self-similar} (ii) and (iii) for $\la=\ga_k$ and $\la=\eta_k$ respectively, which can be essentially deduced from \cite[Theorem 2]{Dajani-Jiang-Kong-Li-Xi-2021}.
%
\begin{proof}[Proof of Theorem \ref{th:unique-coding-totally-self-similar} (ii) and (iii)]
 For (ii) let $\la=\ga_k=\rho(1-\rho^k)$ with $k\in\N$. Then, in view of Figure \ref{Fig:1} (left), we have $f_\la(I)\cap f_{1-\rho-\la}(I)=\emptyset$, where $I=[0,1]$. Furthermore,
  \begin{align*}
 f_0(I)\cap f_\la(I)=[\la, \rho]=[\rho(1-\rho^k), \rho]=f_{\la 0^k}(I)=f_{0(1-\rho)^k}(I),
  \end{align*}
  and symmetrically,
  \[
  f_{1-\rho-\la}(I)\cap f_{1-\rho}(I)=[1-\rho,1-\la]=f_{(1-\rho-\la)(1-\rho)^k}(I)=f_{(1-\rho)0^k}(I).
  \]
  Thus, the IFS $(E_\la, \mathcal F_\la)$ belongs to the class $\mathcal E$ studied in \cite{Dajani-Jiang-Kong-Li-Xi-2021}. So, by \cite[Theorem 2]{Dajani-Jiang-Kong-Li-Xi-2021} it follows that $\dim_H E_\la=t\in(0,1)$ satisfies
  \[
  \rho^t(1-\rho^{kt})+\rho^t+\rho^t(1-\rho^{kt})+\rho^t=1,
  \]
  which can be simplified as $4\rho^t-2\rho^{(k+1)t}=1$. Furthermore, $\dim_H U_\la=s\in(0,1)$ satisfies
  \[
  2\rho^s\left(1-\frac{\rho^{ks}(2-\rho^{ks}-\rho^{ks})}{1-\rho^{2ks}}\right)+2\rho^s=1,
  \]
  which can be deduced as $4\rho^s- \rho^{ks} =1$. This establishes (ii).

  Next we prove (iii). Let $\la=\eta_k=\frac{1-2\rho+\rho^{k+1}}{2}$ for some $k\in\N$. Then, in view of Figure \ref{Fig:1} (right), we have $f_0(I)\cap f_\la(I)=\emptyset$ and $f_{1-\rho-\la}(I)\cap f_{1-\rho}(I)=\emptyset$. Furthermore,
  \[
 f_\la(I)\cap f_{1-\rho-\la}(I)=[1-\rho-\la, \rho+\la]=\left[\frac{1-\rho^{k+1}}{2}, \frac{1+\rho^{k+1}}{2}\right]=f_{\la(1-\rho)^k}(I)=f_{(1-\rho-\la)0^k}(I).
 \]
  So, $(E_\la, \mathcal F_\la)$ also belongs to the class $\mathcal E$ in \cite{Dajani-Jiang-Kong-Li-Xi-2021}. By \cite[Theorem 2]{Dajani-Jiang-Kong-Li-Xi-2021} it follows that $\dim_H E_\la=t\in(0,1)$ satisfies
  \[
  \rho^t+\rho^t(1-\rho^{kt})+\rho^t+\rho^t=1,
  \]
  which can be deduced as $4\rho^t-\rho^{(k+1)t}=1$. Moreover, $\dim_H U_\la=s\in(0,1)$ satisfies
 \[
  \rho^s+\rho^s(1-2\rho^{ks})+\rho^s+\rho^s=1,
  \]
  which can be simplified as $4\rho^s-2\rho^{(k+1)s}=1$. This proves (iii).

  By (i) and (ii) it is easy to verify that $\dim_H U_\la<\dim_H E_\la$, completing the proof.
\end{proof}
\begin{remark}
  \label{rem:41}
  By the above proof it follows that $E_\la$ has an exact overlap when $\la\in\bigcup_{k=1}^\f\set{\ga_k, \eta_k}$.
\end{remark}
Note that when $\la\in\bigcup_{k=1}^\f\set{\ga_k,\eta_k}$, the IFS $(E_\la, \mathcal F_\la)$ belongs to the class $\mathcal E$ studied in \cite{Dajani-Jiang-Kong-Li-Xi-2021}. Then   by \cite[Theorem 1]{Dajani-Jiang-Kong-Li-Xi-2021} it follows that if $\la=\ga_k$ for some $k\in\N$, then  $U_\la$ is not closed, and there are infinitely many $x\in E_\la$ having countably infinitely many codings. On the other hand,  if $\la=\eta_k$ for some $k\in\N$, then $U_\la$ is closed, and there is no $x\in E_\la$ having countably infinitely many codings. Furthermore, by \cite[Theorem 2]{Dajani-Jiang-Kong-Li-Xi-2021} it follows that for any $\la\in\bigcup_{k=1}^\f\set{\ga_k,\eta_k}$ we have
\[\mathcal H^{t}(E_\la)\in(0,+\f)\quad \textrm{and}\quad \mathcal H^{s}(U_\la)\in(0,+\f),\]
where $t=\dim_H E_\la$ and $s=\dim_H U_\la$.

In the following we only need to consider $\la=\la_k=\frac{\rho(1-\rho^k)}{1+\rho^k}$ for some $k\in\N$.
The overlapping  structure is completely different from that for $\la\in\bigcup_{k=1}^\f\set{\ga_k,\eta_k}$ (see Figure \ref{Fig:2}).
In particular, the IFS $(E_\la, \mathcal F_\la)$ does not belong to the class $\mathcal E$ studied  in \cite{Dajani-Jiang-Kong-Li-Xi-2021}. Let $\la=\la_k$ for some $k\in\N$, and set
 \begin{equation}\label{eq:ji-1234}
  \begin{split}
 \i_1&=0(1-\rho)^{k-1}(1-\rho-\lambda ), \hspace{3.9cm}\j_1 =\lambda 0^k;\\
 \i_2&=0(1-\rho)^{k}, \hspace{6.15cm} \j_2 =\lambda 0^{k-1}\lambda;\\
 \i_3&=(1-\rho-\lambda)(1-\rho)^{k}, \hspace{4.5cm}\j_3 =(1-\rho)0^{k-1}\lambda;\\
 \i_4&=(1-\rho-\lambda)(1-\rho)^{k-1}(1-\rho-\lambda),\hspace{2.3cm} \j_4 =(1-\rho)0^{k}.
  \end{split}
 \end{equation}
 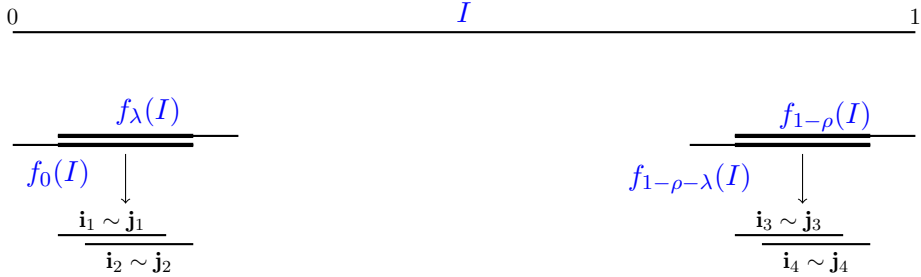
\begin{figure}[h!]
 \begin{center}
 \begin{tikzpicture}[
 scale=6,
 axis/.style={very thick, ->},
 important line/.style={thick},
 dashed line/.style={dashed, thin},
 pile/.style={thick, ->, >=stealth', shorten <=2pt, shorten
 >=2pt},
 every node/.style={color=black}
 ]
 \draw[important line] (0.2, -0.05)--(2.2, -0.05);\node[blue,above,scale=1pt] at(1.2, -0.05){$I$};
 \node[above,scale=0.8pt] at(0.2, -0.05){$0$};\node[above,scale=0.8pt] at(2.2, -0.05){$1$};

 \draw[important line] (0.2, -0.3)--(0.6, -0.3); \draw[important line] (0.3, -0.28)--(0.7, -0.28);
 \node[above,scale=0.8pt] at(0.3, -0.28){};
 \node[blue,above,scale=1pt] at(0.5, -0.28){$f_{\lambda }(I)$};\node[blue,below,scale=1pt] at(0.3, -0.31){$f_0(I)$};
 \draw[ultra thick] (0.3,-0.3)--(0.6,-0.3); \draw[ultra thick] (0.3,-0.28)--(0.6,-0.28);

 \draw[important line] (1.7, -0.3)--(2.1, -0.3); \draw[important line] (1.8, -0.28)--(2.2, -0.28);
 \node[above,scale=0.8pt] at(1.7, -0.39){ };\node[below,scale=0.8pt] at(1.8, -0.2){ };
 \node[blue,above,scale=1pt] at(1.7, -0.43){$f_{1-\rho-\lambda}(I)$};\node[blue,below,scale=1pt] at(2, -0.18){$f_{1-\rho}(I)$};
 \draw[ultra thick] (1.8,-0.3)--(2.1,-0.3); \draw[ultra thick] (1.8,-0.28)--(2.1,-0.28);

 \draw[->](0.45, -0.32)--(0.45, -0.43);
 \draw[important line] (0.3, -0.5)--(0.54, -0.5); \draw[important line] (0.36, -0.52)--(0.6, -0.52);
 \node[above,scale=0.8pt] at(0.42, -0.51){$\i_1\sim\j_1$};\node[below,scale=0.8pt] at(0.48, -0.52){$\i_2\sim\j_2$};

 \draw[->](1.95, -0.32)--(1.95, -0.43);
 \draw[important line] (1.8, -0.5)--(2.04, -0.5); \draw[important line] (1.86, -0.52)--(2.1, -0.52);
\node[above,scale=0.8pt] at(1.92, -0.51){$\i_3\sim\j_3$};\node[below,scale=0.8pt] at(1.98, -0.52){$\i_4\sim\j_4$};

  \end{tikzpicture}
 \end{center}
 \caption{The overlapping structure of $E_\lambda$ with $\lambda=\lambda_k$ for some $k\in\mathbb{N}$.  Then $f_0(E_\la)\cap f_\la(E_\la)=f_{\i_1}(E_\la)\cup f_{\i_2}(E_\la)$ and $f_{1-\rho-\la}(E_\la)\cap f_{1-\rho}(E_\la)=f_{\i_3}(E_\la)\cup f_{\i_4}(E_\la)$ with $f_{\i_\ell}=f_{\j_\ell}$ for all $\ell\in\set{1,2,3,4}$.}\label{Fig:2}
 \end{figure}
 \begin{lemma}\label{lem:f0-cap-flambda}
  Let $\lambda=\lambda_k$ for some $k\in \mathbb{N}$. Then
  {\[
 f_{\i_\ell}=f_{\j_\ell}\quad \textrm{for all} \quad\ell\in\set{1,2,3,4}.
 \]}
  Furthermore,
  \begin{align*}
  f_0(E_{\lambda})\cap f_{\lambda}(E_{\lambda})&=f_{\i_1}(E_{\lambda})\cup f_{\i_2}(E_{\lambda}), \quad
  f_{1-\rho-\la}(E_{\lambda})\cap f_{1-\rho}(E_{\lambda})=f_{\i_3}(E_{\lambda})\cup f_{\i_4}(E_{\lambda}).
  \end{align*}
\end{lemma}
Before proving Lemma \ref{lem:f0-cap-flambda} we point out that the unions in the lemma are NOT disjoint: $f_{\i_1}(E_\la)\cap f_{\i_2}(E_\la)\ne\emptyset$ and $f_{\i_3}(E_\la)\cap f_{\i_4}(E_\la)\ne\emptyset$.
  \begin{proof}
 Note by the symmetry that $E_\lambda$ has the same structure as $E_{1-\rho-\lambda}$. Then it suffices to prove
 \[f_{\i_1}=f_{\j_1},\quad f_{\i_2}=f_{\j_2}\quad \textrm{and}\quad f_0(E_\la)\cap f_\la(E_\la)=f_{\i_1}(E_\la)\cup f_{\i_2}(E_\la).\]
 Since $\la=\lambda_k=\frac{\rho(1-\rho^k)}{1+\rho^k}$,
   by (\ref{eq:ji-1234}) it follows that
    \begin{equation}\label{eq:lam-k-1}
 f_{\i_1}(0)=\rho(1-\rho^{k})-\rho^{k}\lambda =\lambda=f_{\j_1}(0),\quad f_{\i_2}(0)=\rho(1-\rho^{k})=\lambda(1+\rho^{k})=f_{\j_2}(0).
 \end{equation}
   This implies that $f_{\i_\ell}= f_{\j_\ell}$ for $\ell=1,2$.
 Furthermore,  by (\ref{eq:lam-k-1}) it follows that
 \[
  f_{\i_1}(I)\cup f_{\i_2}(I)=[\la, \rho-\rho^k\la]\cup[\rho-\rho^{k+1},\rho]=[\la,\rho]=f_0(I)\cap f_\la(I),
  \]
where the second equality follows by $\la<\rho$.
 Then by  {Definition \ref{def:totoally-self-similar} and Theorem \ref{th:totally-self-similar}} we obtain
 \begin{equation*}
  \begin{split}
 f_{\i_1}(E_{\lambda})\cup f_{\i_2}(E_{\lambda}) &=(f_{\i_1}(I)\cap E_\la)\cup(f_{\i_2}(I)\cap E_\la) \\
  &=(f_{\i_1}(I)\cup f_{\i_2}(I))\cap E_\la \\
  &=(f_0(I)\cap f_\la(I))\cap E_\la=f_0(E_{\lambda})\cap f_{\lambda}(E_{\lambda})
  \end{split}
 \end{equation*}
as desired.
  \end{proof}

 Our next result shows that for $\la=\la_k$, $U_\la$ can be represented as a strongly connected graph-directed set  satisfying the  {SSC}.
 Let $X_\la\subset\Om_\la^\N$ be a subshift of finite type with the set of forbidden blocks given by
\[
  \mathbf F=\bigcup_{\ell=1}^4\set{\i_\ell, \j_\ell},
  \]
 where $\i_\ell, \j_\ell$ are defined in (\ref{eq:ji-1234}). Let $\si$ be the left-shift map on $\Om_\la^\N$.
 \begin{lemma}
   \label{lem:transitive-U}
   Let $\la=\la_k$ for some $k\in\N$. Then
  $(X_\la, \si)$ is a transitive subshift of finite type.
 \end{lemma}
 \begin{proof}
   Note that each forbidden block in $\mathbf F$ has length $k+1$. Then $(X_\la, \si)$ is a $k$-step subshift of finite type.
   By \cite[Theorem 2.1.8]{Lind_Marcus_1995} it suffices to prove that for any two admissible words {$\c=c_1\ldots c_k, \d=d_1\ldots d_k\in B_*(X_\la)$}, we can find a word $w$ such that $\mathbf{c} w  \mathbf{d}\in B_*(X_\la)$.
   Here $B_*(X_\la)$ denotes the set of  all admissible  words appearing in some sequence of $X_\la$. Take $\c=c_1\ldots c_k, \d=d_1\ldots d_k\in B_*(X_\la)$. We will prove in the following two cases  the existence of   $w$ so that $\c w\d\in B_*(X_\la)$.

 {Case I. $c_k\in\set{0,1-\rho-\la}$. If $d_1\in\set{0,\la}$, then by taking $w=(1-\rho-\la)^{k-1}\la(1-\rho-\la)$ one can verify that the longer word $\c w\d$ does not contain any block from $\mathbf F$, i.e., $\c w\d\in B_*(X_\la)$. If $d_1\in\set{1-\rho-\la, 1-\rho}$, then by taking $w=(1-\rho-\la)^{k-1}\la(1-\rho-\la)\la$ we have $\c w\d\in B_*(X_\la)$.}

{Case II. $c_k\in\set{\la,1-\rho}$. If $d_1\in\set{0,\la}$, then by taking $w=\la^{k-1}(1-\rho-\la)\lambda(1-\rho-\la)$ we have $\c w\d\in B_*(X_\la)$. If $d_1\in\set{1-\rho-\la, 1-\rho}$, then by taking $w=\la^{k-1}(1-\rho-\la)\lambda$ one can verify that $\c w\d\in B_*(X_\la)$.}
\end{proof}

\begin{lemma}
  \label{lem:U-lambda-X-lambda}
  Let $\la=\la_k$ for some $k\in\N$. Then $U_\la=\pi_\la( X_\la)$.
\end{lemma}
\begin{proof}
 Take $\la=\la_k$ and let $x\in E_\la\setminus\pi_\la( X_\la)$. Note that $E_\la=\pi_\la(\Om_\la^\N)$. Then $x$ has a coding   $(x_i)\in\Om_\la^\N$ which contains a block from $\mathcal F=\bigcup_{s=1}^4\set{\i_\ell, \j_\ell}$. So by Lemma \ref{lem:f0-cap-flambda} it follows that $x$ has at least two different codings with the substitution: $\i_\ell\sim \j_\ell$ for $\ell\in\set{1,2,3,4}$. Thus, $x\notin U_\la$.

 On the other hand, take $x\in E_\la\setminus U_\la$. Then $x$ has two different codings, say $(c_i), (d_i)\in\Om_\la^\N$. Without loss of generality we may assume $c_1<d_1$. By the overlapping structure of $E_\la$ we have $x\in [\la, \rho]\cup[1-\rho, 1-\la]$, and by symmetry we may assume $x\in[\la,\rho]$. Then $c_1=0$ and $d_1=\la$.
Note that $\sum_{i=1}^{\infty} \rho ^{i-1}c_i=x=\sum_{i=1}^{\infty} \rho ^{i-1}d_i$.
Then
  \begin{equation}\label{uk}
  \sum_{i=2}^{\infty}  \rho ^{i-1}c_i=\sum_{i=2}^{\infty}  \rho ^{i-1}d_i+\lambda.
  \end{equation}

{\bf Claim.}  If $k\ge 2$, then $c_{2}\ldots c_{k}=(1-\rho)^{k-1}$ and $d_{2}\ldots d_{k}=0^{k-1}$.

  Suppose the claim does not hold, and let $\tau\in\set{2,3,\ldots, k}$ be the smallest integer in which $c_\tau\ne 1-\rho$ or $d_\tau\ne 0$.  Then $c_2\ldots c_{\tau-1}=(1-\rho)^{\tau-2}$ and $d_2\ldots d_{\tau-1}=0^{\tau-2}$.
 So, by (\ref{uk}) and $\la=\la_k=\frac{\rho(1-\rho^k)}{1+\rho^k}$ it follows that
\begin{equation}\label{i-s-1}
  \begin{split}
 \sum_{i=\tau}^{\infty} \rho ^{i-\tau}c_{i}-\sum_{i=\tau}^{\infty} \rho ^{i-\tau}d_{i}&= \frac{\lambda}{\rho^{\tau-1}}-\sum_{i=1}^{\tau-2}\frac{1-\rho}{\rho^{i}}=\frac{1-\rho^k}{\rho^{\tau-2}(1+\rho^k)}-\frac{1-\rho^{\tau-2}}{\rho^{\tau-2}}\\
 &=1-\frac{2\rho^{k+2-\tau}}{1+\rho^k}\geq1-\frac{2\rho^2}{1+\rho^k}.
  \end{split}
\end{equation}
Since $\lambda=\frac{\rho(1-\rho^k)}{1+\rho^k}>\frac{2\rho^2}{1+\rho^k}$, by (\ref{i-s-1}) we have
 \begin{equation*}\label{k>11-rho-lambda-0}
 \sum_{i=\tau}^{\infty} \rho ^{i-\tau}c_{i}-\sum_{i=\tau}^{\infty} \rho ^{i-\tau}d_{i}>1-\la= f_{1-\rho}(1)-f_{\lambda}(0),
  \end{equation*}
which implies that $c_{\tau}=1-\rho$ and $d_{\tau}=0$.
This leads to a contradiction with  the definition of $\tau$, and thus proves the claim.

  By (\ref{uk}) and the claim it follows that for $k\in\N$,
\begin{equation*}
  \sum_{i=k+1}^{\infty} \rho ^{i-k-1}c_{i}-\sum_{i=k+1}^{\infty} \rho ^{i-k-1}d_{i}=\frac{\lambda}{\rho^{k}}-\sum_{i=1}^{k-1} \frac{1-\rho}{\rho ^{i}}=1-\frac{2\rho}{1+\rho^k}>\rho+\lambda= f_{\lambda}(1)-f_{0}(0),
  \end{equation*}
  where the inequality holds by using $\la=\frac{\rho(1-\rho^k)}{1+\rho^k}$ and $0<\rho\le 1/4$.
This implies that $c_{k+1}\in \left\{1-\rho-\lambda,1-\rho\right\}$ and $d_{k+1}\in \left\{0,\lambda\right\}$.
Thus, by (\ref{eq:ji-1234}) it follows that $c_1 c_2\cdots c_{k+1}\in \left\{\i_1,\i_2\right\}$ and $d_1 d_2\cdots d_{k+1}\in \left\{\j_1,\j_2\right\}$, which yields
\[(c_i), (d_i)\notin X_\la.\]
Since $(c_i), (d_i)$ are two arbitrary codings of $x$, we conclude that $x\notin\pi_\la(X_\la)$.
\end{proof}

\begin{lemma}
  \label{lem:unique-coding-lambdak-1}
  Let $\la=\la_k$ for some $k\in\N$. Then $U_\la$ can be represented as a strongly connected graph-directed set satisfying the {SSC}.
\end{lemma}
\begin{proof}
  Note that $X_\la$ is a subshift of finite type with the set $\mathbf F=\bigcup_{\ell=1}^4\set{\i_\ell, \j_\ell}$ of forbidden blocks.
   {Since each block in $\mathbf F$ has length $k+1$, $X_\la$ is a $k$-step subshift of finite type which can be represented as a directed graph $\mathcal G=(\mathcal V, \mathcal E)$ constructed in the following way.}
  Let $\mathcal V=B_k(X_\la)$ be the set of all length $k$ admissible blocks appearing in some sequence of $X_\la$.
  For two vertices $\c=c_1\ldots c_k, \d=d_1\ldots d_k\in\mathcal V$ we draw a directed edge from $\c$ to $\d$, denoted by $\overrightarrow{\c\d}$, if
  \[
  c_2\ldots c_k=d_1\ldots d_{k-1}\quad\textrm{and}\quad c_1\ldots c_kd_k\in B_{k+1}(X_\la).
  \]
In this case we define a map for this edge $\overrightarrow{\c\d}$ by $f_{\overrightarrow{\c\d}}=f_{c_1}$.
Let $\mathcal E$ be the set of all directed edges in $\mathcal G$.

  For $\c=c_1\ldots c_k\in\mathcal V$ we set
  \[
  U_\c:=\set{\pi_\la((x_i)): (x_i)\in X_\la\quad\textrm{and}\quad x_1\ldots x_k=c_1\ldots c_k}.
  \]
  Then by Lemma  \ref{lem:U-lambda-X-lambda} it follows that
  \[
  U_\la=\pi_\la(X_\la)=\bigcup_{\c\in\mathcal V}U_\c,
  \]
  where
  \begin{equation}\label{eq:U-ssc-1}
  U_\c=\bigcup_{\overrightarrow{\c\d}\in\mathcal E}f_{\overrightarrow{\c\d}}(U_\d).
  \end{equation}
  Note by Lemma \ref{lem:transitive-U} that the graph $\mathcal G=(\mathcal V, \mathcal E)$ is strongly connected. Then it suffices to prove that the union in (\ref{eq:U-ssc-1}) is pairwise disjoint.

  Suppose on the contrary there exist $\d=d_1\ldots d_k, \mathbf e=e_1\ldots e_k\in\mathcal V$ such that $\overrightarrow{\c\d}, \overrightarrow{\c\mathbf e}\in\mathcal E$ and $f_{\overrightarrow{\c\d}}(U_\d)\cap f_{\overrightarrow{\c\mathbf e}}(U_{\mathbf e})\ne\emptyset$. Then by the definition of $\mathcal G=(\mathcal V, \mathcal E)$ it follows that
  \begin{equation}\label{eq:U-ssc-2}
  d_1\ldots d_{k-1}=c_2\ldots c_k=e_1\ldots e_{k-1}\quad\textrm{and}\quad d_k\ne e_k.
  \end{equation}
  Furthermore, there exist $x=\pi_\la((x_i))\in U_\d, y=\pi_\la((y_i))\in U_{\mathbf e}$ such that $f_{c_1}(x)=f_{c_1}(y)$.
  Note by Lemma \ref{lem:U-lambda-X-lambda} that  {$U_{\mathbf d}, U_{\mathbf e}\subset\pi_\la(X_\la)=U_\la$}. Then $(x_i), (y_i)$ are the unique codings of $x$ and $y$, respectively. Furthermore,  $x_1\ldots x_{k}=d_1\ldots d_k$ and $y_1\ldots y_k=e_1\ldots e_k$. By (\ref{eq:U-ssc-2}) and using $f_{c_1}(x)=f_{c_1}(y)$ we obtain that
  \[
 z:= \sum_{i=1}^\f\rho^{i-1}x_{k-1+i}=\sum_{i=1}^\f\rho^{i-1}y_{k-1+i}
  \]
  has   two different codings $x_k x_{k+1}\ldots$ and $y_ky_{k+1}\ldots$. This leads to a contradiction with Lemma \ref{lem:U-lambda-X-lambda} that $z=\pi_\la(x_kx_{k+1}\ldots)\in\pi_\la(X_\la)=U_\la$ should have a unique coding.
\end{proof}
\begin{proposition}
  \label{prop:unique-coding-lambdak}
  Let $\la=\la_k$ for some $k\in\N$. Then
\[
  \dim_H U_\la=s,
  \]
 where $s\in(0,1)$ satisfies $4\rho^s-2\rho^{ks}=1$.
\end{proposition}
\begin{proof}
Take $\la=\la_k$.
For a word $w=w_1\ldots w_n\in B_n(X_\la)$ let $U_\la(w):=U_\la\cap f_{w}(E_\la)$. Then $U_\la(w)$ consists of all $x\in U_\la$ whose unique coding beginning with the word $w$. So,
\begin{equation}\label{eq:janu-18-1}
   U_\la=\bigcup_{d\in\Om_\la}U_\la(d)
\end{equation}
with the union pairwise disjoint.
Note by Lemma \ref{lem:unique-coding-lambdak-1}  that $U_\la$ is a strongly connected graph-directed set satisfying the SSC. Then by \cite{Mauldin_Williams_1988} it follows that for $s=\dim_H U_\la$ we have $\mathcal H^s(U_\la)\in(0,\f)$. So by (\ref{eq:janu-18-1}) it suffices to prove that for each $d\in\Om_\la=\set{0,\la, 1-\rho-\la, 1-\rho}$,
\begin{equation}
  \label{eq:janu-18-2}
  \mathcal H^s(U_\la(d))=\left(\rho^s-\frac{2\rho^{(k+1)s}}{1+2\rho^{ks}}\right)\mathcal H^s(U_\la).
\end{equation}
Since the proofs of (\ref{eq:janu-18-2}) for different $d\in\Om_\la$ are similar, we only prove (\ref{eq:janu-18-2}) for $d=0$.

Note by Lemma \ref{lem:U-lambda-X-lambda} that $U_\la=\pi_\la(X_\la)$.
This means that for a given   $x\in E_\la$, if all of its codings do not contain any block from $\mathbf F=\bigcup_{\ell=1}^4\set{\i_\ell,\j_\ell}$ then $x\in U_\la$. So by the definition of $U_\la(w)$ we obtain that
\begin{equation}
  \label{eq:janu-18-3}
  U_\la(0)=f_0(U_\la)\setminus\bigcup_{\ell=1}^2 f_0(U_\la(\hat\i_\ell)),
\end{equation}
where for $\ell=1,2,3,4$ we set $\hat\i_\ell:=\si(\i_\ell)$ and $\hat\j_\ell:=\si(\j_\ell)$. Note by (\ref{eq:ji-1234}) that $\hat\i_{5-\ell}=\hat\i_\ell$ and $\hat\j_{5-\ell}=\hat\j_\ell$ for any $\ell\in\set{1,2,3,4}$. Then by the definition of $U_\la(w)$ and using $U_\la=\pi_\la(X_\la)$ it follows that
\begin{equation}\label{eq:janu-18-4}
\begin{split}
  U_\la(\hat\i_1)&=f_{\hat\i_1}(U_\la)\setminus\bigcup_{\ell=1}^2 f_{\hat\i_1}(U_\la(\hat\i_\ell)),\\
U_\la(\hat\i_2)&=f_{\hat\i_2}(U_\la)\setminus\bigcup_{\ell=1}^2 f_{\hat\i_2}(U_\la(\hat\j_\ell)),\\
  U_\la(\hat\j_1)&= f_{\hat\j_1}(U_\la)\setminus\bigcup_{\ell=1}^2 f_{\hat\j_1}(U_\la(\hat\i_\ell)),\\
  U_\la(\hat\j_2)&=f_{\hat\j_2}(U_\la)\setminus\bigcup_{\ell=1}^2 f_{\hat\j_2}(U_\la(\hat\j_\ell)).
\end{split}
\end{equation}
Since the unions in (\ref{eq:janu-18-4}) are pairwise disjoint, by taking the $s$-dimensional Hausdorff measure on both sides of (\ref{eq:janu-18-4}) we obtain that
\[
\mathcal H^s(U_\la(\hat\i_\ell))=\mathcal H^s(U_\la(\hat\j_\ell))=\frac{\rho^{ks}}{1+2\rho^{ks}}\mathcal H^s(U_\la)\quad\forall ~\ell=1,2.
\]
Thus, by (\ref{eq:janu-18-3}) it follows that
\[
\mathcal H^s(U_\la(0))=\rho^s\left[\mathcal H^s(U_\la)-\mathcal H^s(U_\la(\hat\i_1))-\mathcal H^s(U_\la(\hat\i_2))\right]=\left(\rho^s-\frac{2\rho^{(k+1)s}}{1+2\rho^{ks}}\right)\mathcal H^s(U_\la)
\]
proving (\ref{eq:janu-18-2}) for $d=0$. This completes the proof.
\end{proof}

\begin{proof}  [Proof of Theorem \ref{th:unique-coding-totally-self-similar} (i)]
By Proposition \ref{prop:unique-coding-lambdak} it suffices to prove that for $\la=\la_k$ we have $\dim_H U_\la<\dim_H E_\la$. Let $\mathbf F'=\mathbf F\setminus\set{\i_1}=\set{\i_2, \i_3,\i_4,\j_1,\j_2,\j_3,\j_4}$, and let $X_\la'$ be the subshift of finite type over $\Om_\la$ with the set $\mathbf F'$ of forbidden blocks. Then
\[
X_\la'=\set{(d_i)\in\Om_\la^\N: d_{i+1}\ldots d_{i+k+1}\notin{\mathbf F'}~\forall i\ge 0}.
\]
By a similar argument as in the proof of Lemmas \ref{lem:transitive-U} and \ref{lem:unique-coding-lambdak-1} one can show that $\pi_\la(X_\la')$ is a strongly connected graph-directed set satisfying the SSC. Note that $X_\la$ is a proper subset of $X_\la'$. Then by \cite[Corollary 4.4.9]{Lind_Marcus_1995} it follows that $h_{top}(X_\la)<h_{top}(X_\la')$. This implies that
\[
\dim_H\pi_\la(X_\la)=\frac{h_{top}(X_\la)}{-\log\rho}<\frac{h_{top}(X_\la')}{-\log\rho}=\dim_H\pi_\la(X_\la').
\]
Since $U_\la=\pi_\la(X_\la)$ by Lemma \ref{lem:U-lambda-X-lambda}   and $\pi_\la(X_\la')\subset E_\la$, we conclude that
\[
\dim_H U_\la<\dim_H E_\la,
\]
completing the proof.
\end{proof}
\section{Typical result for the Hausdorff dimension of $U_\la$}\label{sec:typical-dimension-U}

When $E_\la$ is totally self-similar, we determine the Hausdorff dimension of $U_\la$ in the previous section. In this section we show that $U_\la$ has the same Hausdorff dimension as $E_\la$ for typical $\la$, and prove Theorem \ref{th:unique-codings-typical}. Note that
\begin{equation}\label{eq:E-U}
E_\la\setminus U_\la=\bigcup_{\i\in\Om_\la^*}f_\i(M_\la),
\end{equation}
where
\[
M_\la:=\bigcup_{ c,d\in\Om_\la, c\ne d}f_c(E_\la)\cap f_d(E_\la).
\]
Then $E_\la\setminus U_\la$ is a countable union of scaling copies of $M_\la$.
{By the countable stability  of Hausdorff dimension, to prove Theorem \ref{th:unique-codings-typical} (i) it suffices to prove that for $\rho\in(0,1/4)$ we have $\dim_H M_\la<\dim_H E_\la$ for Lebesgue almost every $\la\in(0,\rho)\cup(\frac{1-2\rho}{2}, \frac{1-\rho}{2})$.}
Moreover, to prove Theorem \ref{th:unique-codings-typical} (ii) we only need to prove that for $\rho\in(0, 1/16)$ we have $M_\la=\emptyset$ for Lebesgue almost every $\la\in(0,\rho)\cup(\frac{1-2\rho}{2}, \frac{1-\rho}{2})$. Since the proof  for $\la\in (\frac{1-2\rho}{2}, \frac{1-\rho}{2})$ is similar, we only prove it for $\la\in(0,\rho)$.

Let $J:=[a,b]\subset(0,\rho)$, and take $\la\in J$. Then (see the left graph of Figure \ref{Fig:1})
\[
M_\la=\left(f_0(E_\la)\cap f_\la(E_\la)\right)\cup\left(f_{1-\rho-\la}(E_\la)\cap f_{1-\rho}(E_\la)\right).
\]
Note by symmetry that $f_{1-\rho-\la}(E_\la)\cap f_{1-\rho}(E_\la)=1-f_0(E_\la)\cap f_\la(E_\la)$. So in the following it suffices to prove that for $\rho\in(0, 1/4)$ we have $\dim_H(f_0(E_\la)\cap f_\la(E_\la))<\dim_H E_\la$ for Lebesgue almost every $\la\in J$; and for $\rho\in(0,1/16)$ we have  $f_0(E_\la)\cap f_\la(E_\la)=\emptyset$ for Lebesgue almost every $\la\in J$. Suppose $f_0(E_\la)\cap f_\la(E_\la)\ne\emptyset$ for some $\la\in J$. Otherwise, we are done. Observe that each $x\in E_\la$ has a coding $\i\in\Om_\la^\N$ satisfying  $x=\pi_\la(\i)$. Set
\[
D_J:=\set{(\i,\j)\in\Om_\la^\N\times\Om_\la^\N: \exists \la\in J\textrm{ such that }f_0(\pi_\la(\i))=f_\la(\pi_\la(\j))}.
\]
\begin{lemma}
  \label{lem:existence-lammda}
  Each pair $(\i,\j)\in D_J$ determines a unique $\la\in J$.
\end{lemma}
\begin{proof}
  Take $(\i,\j)\in D_J$ with $\i=i_1i_2\ldots, \j=j_1j_2\ldots$. Then $f_0(\pi_\la(\i))=f_\la(\pi_\la(\j))$. By (\ref{eq:coding-map}) it follows that
  \begin{equation}
 \label{eq:unique-1}
 \sum_{n=1}^\f\rho^n i_n=\sum_{n=1}^\f\rho^n j_n+\la.
  \end{equation}
 Note that the digits $i_n, j_n\in\Om_\la=\set{0,\la,1-\rho-\la,1-\rho}$ might contain the parameter $\la$. In order to separate the parameter $\la$ we partition the set $\N$ into
$
 \nn_\i^1:=\set{n: i_n=0}, \nn_\i^2:=\set{n: i_n=\la}, \nn_\i^3:=\set{n: i_n=1-\rho-\la}$ and $\nn_\i^4:=\set{n: i_n=1-\rho}.
$ Similarly, for $s\in\set{1,2,3,4}$ we define $\nn_\j^s$ by replacing the $i_n$ in $\nn_\i^s$ by $j_n$. Thus, (\ref{eq:unique-1}) can be rewritten as
\[
\la\sum_{n\in\nn_\i^2}\rho^n+(1-\rho-\la)\sum_{n\in\nn_\i^3}\rho^n+(1-\rho)\sum_{n\in\nn_\i^4}\rho^n= \la\sum_{n\in\nn_\j^2}\rho^n+(1-\rho-\la)\sum_{n\in\nn_\j^3}\rho^n+(1-\rho)\sum_{n\in\nn_\j^4}\rho^n+\la,
\]
which can be reorganized as
\[
\la\left(1+\sum_{n\in\nn_\j^2}\rho^n+\sum_{n\in\nn_\i^3}\rho^n-\sum_{n\in\nn_\j^3}\rho^n-\sum_{n\in\nn_\i^2}\rho^n\right)=(1-\rho)\left(\sum_{n\in\nn_\i^3\cup\nn_\i^4}\rho^n-\sum_{n\in\nn_\j^3\cup\nn_\j^4}\rho^n\right).
\]
Since $0<\rho<1/4$, we have
\begin{equation}
  \label{eq:unique-2}
  1+\sum_{n\in\nn_\j^2}\rho^n+\sum_{n\in\nn_\i^3}\rho^n-\sum_{n\in\nn_\j^3}\rho^n-\sum_{n\in\nn_\i^2}\rho^n\ge1-2\sum_{n=1}^\f\rho^n=\frac{1-3\rho}{1-\rho}>0.
\end{equation}
This gives
\[
\la=\frac{(1-\rho)\left(\sum_{n\in\nn_\i^3\cup\nn_\i^4}\rho^n-\sum_{n\in\nn_\j^3\cup\nn_\j^4}\rho^n\right)}{1+\sum_{n\in\nn_\j^2}\rho^n+\sum_{n\in\nn_\i^3}\rho^n-\sum_{n\in\nn_\j^3}\rho^n-\sum_{n\in\nn_\i^2}\rho^n}
\]
as desired.
\end{proof}

In terms of Lemma \ref{lem:existence-lammda}, let $\la_{\i,\j}$ be the unique $\la\in J$ determined by the pair $(\i,\j)\in D_J$. Then
\begin{equation}\label{eq:lambda-ij}
\la_{\i,\j}=\frac{(1-\rho)p_{\i,\j}}{1+q_{\i,\j}},
\end{equation}
where
\[
p_{\i,\j}:=\sum_{n\in\nn_\i^3\cup\nn_\i^4}\rho^n-\sum_{n\in\nn_\j^3\cup\nn_\j^4}\rho^n,\quad q_{\i,\j}:=\sum_{n\in\nn_\j^2}\rho^n+\sum_{n\in\nn_\i^3}\rho^n-\sum_{n\in\nn_\j^3}\rho^n-\sum_{n\in\nn_\i^2}\rho^n.
\]
Equipped with the metric $d$ on $\Om_\la^\N$ given by
\[
d(\i,\j):=\rho^{\inf\set{n: i_n\ne j_n}},
\]
we define a metric $\|\cdot\|$ on the product space $\Om_\la^\N\times\Om_\la^\N$ by
\[
\|(\i,\j), (\u,\v)\|=\max\set{d(\i,\u), d(\j,\v)}.
\]
\begin{lemma}
  \label{lem:lipschitz}
  The map $\Phi: D_J\to J\times [0,1]$ defined by
  \[
  \Phi((\i,\j))=(\la_{\i,\j}, f_0(\pi_{\la_{\i,\j}}(\i)))
  \]
  is Lipschitz continuous with respect to the metric $\|\cdot\|$ on $D_J$.
\end{lemma}
\begin{proof}
  Take two pairs $(\i,\j), (\u,\v)\in D_J$. It suffices to prove that
  \begin{equation}
 \label{eq:lip-1}
 |\la_{\i,\j}-\la_{\u,\v}|\le C_1\|(\i,\j), (\u,\v)\|
  \end{equation}
  and
  \begin{equation}
 \label{eq:lip-2}
 |f_0(\pi_{\la_{\i,\j}}(\i))-f_0(\pi_{\la_{\u,\v}}(\u))|\le C_2\|(\i,\j), (\u,\v)\|
  \end{equation}
  for some constants $C_1, C_2>0$. Note by (\ref{eq:lambda-ij}) that
  \begin{equation}
 \label{eq:lip-3}
 \begin{split}
 \left|\frac{\la_{\i,\j}-\la_{\u,\v}}{1-\rho}\right|&=\left|\frac{p_{\i,\j}}{1+q_{\i,\j}}-\frac{p_{\u,\v}}{1+q_{\u,\v}}\right|\\
 &\le\frac{1}{|1+q_{\i,\j}|}|p_{\i,\j}-p_{\u,\v}|+\left|\frac{p_{\u,\v}}{(1+q_{\i,\j})(1+q_{\u,\v})}\right|\cdot|q_{\i,\j}-q_{\u,\v}|\\
 &\le\frac{1-\rho}{1-3\rho}|p_{\i,\j}-p_{\u,\v}|+\frac{\rho(1-\rho)}{(1-3\rho)^2}|q_{\i,\j}-q_{\u,\v}|,
 \end{split} \end{equation}
 where the last inequality follows by (\ref{eq:unique-2}) that $1+q_{\i,\j}, 1+q_{\u,\v}\ge \frac{1-3\rho}{1-\rho}>0$ and $|p_{\u,\v}|\le\sum_{n=1}^\f\rho^n=\frac{\rho}{1-\rho}$.
  Observe that
  \begin{align*}
 |p_{\i,\j}-p_{\u,\v}|&\le\left|\sum_{n\in\nn_\i^3\cup\nn_\i^4}\rho^n-\sum_{n\in\nn_\u^3\cup\nn_\u^4}\rho^n\right|+\left|\sum_{n\in\nn_\j^3\cup\nn_\j^4}\rho^n-\sum_{n\in\nn_\v^3\cup\nn_\v^4}\rho^n\right|\\
 &\le \tilde C_1[d(\i,\u)+d(\j,\v)]\le \tilde C_2\|(\i,\j), (\u,\v)\|
  \end{align*}
  for some constants $\tilde C_1, \tilde C_2>0$, and similarly, $|q_{\i,\j}-q_{\u,\v}|\le \tilde C_3\|(\i,\j), (\u,\v)\|$ for some constant $\tilde C_3>0$. So, by (\ref{eq:lip-3}) we prove (\ref{eq:lip-1}).

  On the other hand, observe by $\la_{\i,\j} \in J=[a,b]$ that
  \begin{align*}
  |f_0(\pi_{\la_{\i,\j}}(\i))-f_0(\pi_{\la_{\u,\v}}(\u))|&\le (1-\rho)\left|\sum_{n\in\nn_\i^3\cup\nn_\i^4}\rho^n-\sum_{n\in\nn_\u^3\cup\nn_\u^4}\rho^n\right|\\
  &\quad+\left|\la_{\i,\j}\left(\sum_{n\in\nn_\i^2}\rho^n-\sum_{n\in\nn_\i^3}\rho^n\right)-\la_{\u,\v}\left(\sum_{n\in\nn_\u^2}\rho^n-\sum_{n\in\nn_\u^3}\rho^n\right)\right|\\
  &\le  \tilde C_4 \cdot d(\i,\u)+|\la_{\i,\j}|\cdot\left|\sum_{n\in\nn_\i^2}\rho^n+\sum_{n\in\nn_\u^3}\rho^n-\sum_{n\in\nn_\i^3}\rho^n-\sum_{n\in\nn_\u^2}\rho^n\right|\\
  &\quad +\left|\sum_{n\in\nn_\u^2}\rho^n-\sum_{n\in\nn_\u^3}\rho^n\right|\cdot|\la_{\i,\j}-\la_{\u,\v}|\\
  &\le \tilde C_4\cdot d(\i,\u)+\tilde C_5\cdot d(\i,\u)+\frac{\rho}{1-\rho}|\la_{\i,\j}-\la_{\u,\v}|\\
  &\le C_2\|(\i,\j), (\u,\v)\|
  \end{align*}
  for some constant $C_2>0$, where the last inequality follows by (\ref{eq:lip-1}). This proves (\ref{eq:lip-2}), completing the proof.
\end{proof}

The  following famous slicing theorem is due to Marstrand \cite{Marstrand-1954}.
\begin{lemma}
  \label{1lem:Marstrand-slicing-theorem}
  Let $K\subset \R^2$ be a Borel set. Then for Lebesgue almost every $x\in\R$ we have
  \[
  \dim_H K\cap\set{(x,y): y\in\R}\le \max\set{0,\dim_H K-1}.
  \]
  In particular, if $\dim_H K<1$, then for Lebesgue almost every $x\in \R$ the intersection
  \[K\cap\set{(x,y): y\in\R}=\emptyset.\]
\end{lemma}

Note that $E_\la$ is an affine image of the projection of the four corner Cantor set $E$ generated by the IFS
 $\set{(\rho x, \rho y), (\rho x, \rho y+1-\rho), (\rho x+1-\rho, \rho y), (\rho x+1-\rho, \rho y+1-\rho)}$.
Another useful result was essentially due to Hochman \cite{Hochman-2014} (see also \cite[Theorem 2.1]{Shmerkin-2015}).
\begin{lemma}
  \label{1lem:dimension}
  Let $\rho\in(0,1/4)$. Then for Lebesgue almost every $\la\in(0,\frac{1-\rho}{2})$ we have
  \[
  \dim_H E_\la=\frac{2\log 2}{-\log \rho}.
  \]
\end{lemma}
 \begin{proof}
   [Proof of Theorem \ref{th:unique-codings-typical}]
  First we consider $\rho\in(0, 1/4)$. Take $J=[a, b]\subset(0,\rho)$. By Lemma \ref{lem:lipschitz} it follows that
  \begin{equation}\label{eq:dim-Phi-DJ}
  \dim_H\Phi(D_J)\le \dim_H D_J\le 2\dim_H \Om_\la^\N=2\frac{\log 4}{-\log \rho},
  \end{equation}
  where the last equality holds since $\Om_\la^\N$ is a compact metric space under the metric $d(\i,\j)=\rho^{\inf\set{n: i_n\ne j_n}}$.
Note that
$
M_\la=[f_0(E_\la)\cap f_\la(E_\la)]\cup[1-f_0(E_\la)\cap f_\la(E_\la)].
$ So, by Lemma \ref{1lem:Marstrand-slicing-theorem} it follows that for Lebesgue almost every $\la\in J$,
\begin{align*}
  \dim_H M_\la&=\dim_H(f_0(E_\la)\cap f_\la(E_\la))\le\dim_H(\Phi(D_J)\cap\set{(\la, y): y\in\R})\\
  &\le \dim_H\Phi(D_J)-1\le \frac{2\log 4}{-\log\rho}-1<\frac{ \log 4}{-\log \rho},
 \end{align*}
 where the last inequality follows by  {$0<\rho<1/4$}.
 Hence, by (\ref{eq:E-U}) and Lemma \ref{1lem:dimension} it follows that
 for Lebesgue almost every $\la\in J$,
 \[
 \dim_H (E_\la\setminus U_\la)<\frac{\log 4}{-\log\rho}=\dim_H E_\la,
 \]
 which yields $\dim_H U_\la=\dim_H E_\la$. Since $J\subset(0,\rho)$ was arbitrary, it follows that $\dim_H U_\la=\dim_H E_\la=\frac{\log 4}{-\log\rho}$ for Lebesgue almost every $\la\in(0,\rho)$.

  Next we assume $\rho\in(0, 1/16)$.   Then by (\ref{eq:dim-Phi-DJ}) we have $\dim_H \Phi(D_J)\le \frac{2\log 4}{-\log \rho}<1$.
  So, by the second statement of Lemma \ref{1lem:Marstrand-slicing-theorem} one can deduce that for Lebesgue almost every $\la\in J$  the intersection $f_0(E_\la)\cap f_\la(E_\la)=\emptyset$, and then  $M_\la=\emptyset$. So, by (\ref{eq:E-U}) we conclude that
  $U_\la=E_\la$ for Lebesgue almost every $\la\in J$. Since $J\subset(0,\rho)$ was arbitrary,   it follows that $U_\la=E_\la$ for Lebesgue almost every $\la\in(0, \rho)$.
 \end{proof}

\section{The possibility for $E_{\la}$ to contain  an interval}\label{sec:density}

Recall by (\ref{eq:W}) that $W$ consists of all coprime pairs $(p, q)\in\N^2$ with $p<q$ and $ord_2(p), ord_2(q)$ even.
 By Theorem \ref{th:martilla} it follows that  for $\la\in(0,3/8)$   the self-similar set $E_{\la}$ has an exact overlap  if and only if $\la=\frac{3p}{4(p+q)}$ with  $(p,q)\in W$. Similarly, we recall from (\ref{eq:hat-W}) that $\hat W$ consists of all coprime pairs $(p, q)\in\N^2$ satisfying $p<q$, $ord_2(p)$ odd or $ord_2(q)$ odd. Furthermore, by Theorem \ref{th:martilla} it follows that for $\la\in(0, 3/8)$, $E_\la$ contains a non-degenerate interval if and only if $\la=\frac{3p}{4(p+q)}$ with $(p,q)\in \hat W$.
 In this section we will
 describe the densities of $W$ and $\hat W$ in $\mathbb{N}^{2}$, and prove Theorem \ref{th:density-exact-overlap}.

 First  we recall some known results from analytic number theory (cf.~\cite{Hardy-Wright-2008}).
Let $\phi$ be the  Euler's  function such that for $n\in\N$,  $\phi(n)$ is the number of positive integers no larger than and prime to $n$.  Then $\phi(1)=1$, and for  $n\in\N_{\ge 2}$, if we write it  in a standard form $n=p^{c_{1}}_{1}p^{c_{2}}_{2}\cdots p^{c_{r}}_{r}$ with $p_1, p_2,\ldots, p_r$   distinct primes, then  (cf.~\cite[Theorem 62]{Hardy-Wright-2008})
\begin{equation}
  \label{Euler's totient function}
  \phi (n)=\prod_{i=1}^{r}p^{c_{i}-1}_{i}(p_{i}-1).
\end{equation}
Furthermore,   the summation $\sum_{n=1}^N\phi(n)$ increases to infinity  of order $N^2$. In fact, by \cite[Theorem 330]{Hardy-Wright-2008} we have
\begin{equation}
  \label{eq:sum-Euler-function}
  \lim_{N\to\f}\frac{1}{N^2}\sum_{n=1}^N \phi(n)=\frac{3}{\pi^2}.
\end{equation}
Another useful representation of $\phi$ is based on
  the M\"obius function $\mu$  defined by
\[
\mu(n)=\left\{\begin{array}
  {lll}
  1&\textrm{if}& n=1\\
  0&\textrm{if}& n\textrm{ has a squared factor}\\
  (-1)^k&\textrm{if}& n\textrm{ is the product of $k$ different primes}.
\end{array}\right.
\]
Then by (\ref{Euler's totient function}) the  function $\phi$ can be rewritten as
\begin{equation}\label{eq:Euler-function-1}
{\phi(n)=n\sum_{m|n}\frac{\mu(m)}{m},}
\end{equation}
where the summation is taken over all positive factors $m$ of $n$.

The following result can be easily deduced from \cite[Theorem 287]{Hardy-Wright-2008}.
\begin{lemma}
  \label{lem:zeta-function}
  For any $s>1$ let $
  \zeta(s)=\sum_{n=1}^\f{n^{-s}}
  $ be the zeta function. Then
  \[
  \sum_{n=1}^\f\frac{\mu(n)}{n^s}=\prod_{p}(1-\frac{1}{p^s})=\frac{1}{\zeta(s)},
  \]
  where the product is taken over all prime numbers.
\end{lemma}
Next we prove a useful lemma which is comparable with (\ref{eq:sum-Euler-function}).
\begin{lemma}
  \label{lem:sum-Euler-function-odd}
  \[\lim_{N\to\f}\frac{1}{N}\sum_{n=1}^N\frac{\phi(2n-1)}{2n-1}=\frac{8}{\pi^2}.\]
\end{lemma}
\begin{proof}
  By (\ref{eq:Euler-function-1}) it follows that
  {\begin{align*}
 \sum_{n=1}^N\frac{\phi(2n-1)}{2n-1}
 &=\sum_{n=1}^N\sum_{m|(2n-1)}\frac{\mu(m)}{m}=\sum_{n=1}^N\sum_{(2m-1)|(2n-1)}\frac{\mu(2m-1)}{2m-1}\\
 &=\sum_{m=1}^N\frac{\mu(2m-1)}{2m-1}\sum_{n=1}^N\mathbb I_{\set{(2n-1)(2m-1)\le 2N-1}}\\
 &=\sum_{m=1}^N\frac{\mu(2m-1)}{2m-1}\flr{\frac{N+m-1}{2m-1}}\\
 &=\sum_{m=1}^N\frac{\mu(2m-1)}{2m-1}\cdot\frac{N}{2m-1}+\ep_N,
  \end{align*}
where $\mathbb I$ is the indicator function and
\[
|\ep_N|\le\sum_{m=1}^N\frac{m-1}{(2m-1)^2}+\sum_{m=1}^N\frac{1}{2m-1}<\sum_{m=1}^N\frac{2}{2m-1}.
\]
Clearly, $\frac{|\ep_N|}{N}\to 0$ as $N\to\f$. This implies that
\begin{equation*}\label{eq:sum-Euler-1}
  \lim_{N\to\f}\frac{1}{N}\sum_{n=1}^N\frac{\phi(2n-1)}{2n-1}
  =\lim_{N\to\f} \sum_{m=1}^N\frac{\mu(2m-1)}{(2m-1)^2}
  = \sum_{m=1}^\f\frac{\mu(2m-1)}{(2m-1)^2}.
\end{equation*}
Therefore, the lemma follows by  Lemma \ref{lem:zeta-function}     that
\begin{align*}
  \sum_{m=1}^\f\frac{\mu(2m-1)}{(2m-1)^2}&=\prod_{p\ge 3}\left(1+\frac{\mu(p)}{p^2}+\frac{\mu(p^2)}{p^4}+\cdots\right) =\prod_{p\ge 3}\left(1-\frac{1}{p^2}\right)\\
  &=\frac{1}{1-\frac{1}{2^2}}\prod_{p}\left(1-\frac{1}{p^2}\right)=\frac{4}{3}\frac{1}{\zeta(2)}=\frac{8}{\pi^2},
\end{align*}
where the first two products are taken over all primes at least three,  the third product is taken over all primes,  and the last equality follows by using $\zeta(2)=\sum_{n=1}^\f{n^{-2}}={\pi^2}/{6}$.}
\end{proof}

\begin{proposition}
  \label{prop:density}
  Let $\hat W$ be defined as in (\ref{eq:hat-W}). Then
    \[
  \lim_{N\to\f}\frac{\#(\hat W\cap[1, N]^2)}{N^2}=\frac{4}{3\pi^2}.
  \]
\end{proposition}\begin{figure}[h!]
  \centering
  \includegraphics[width=10cm]{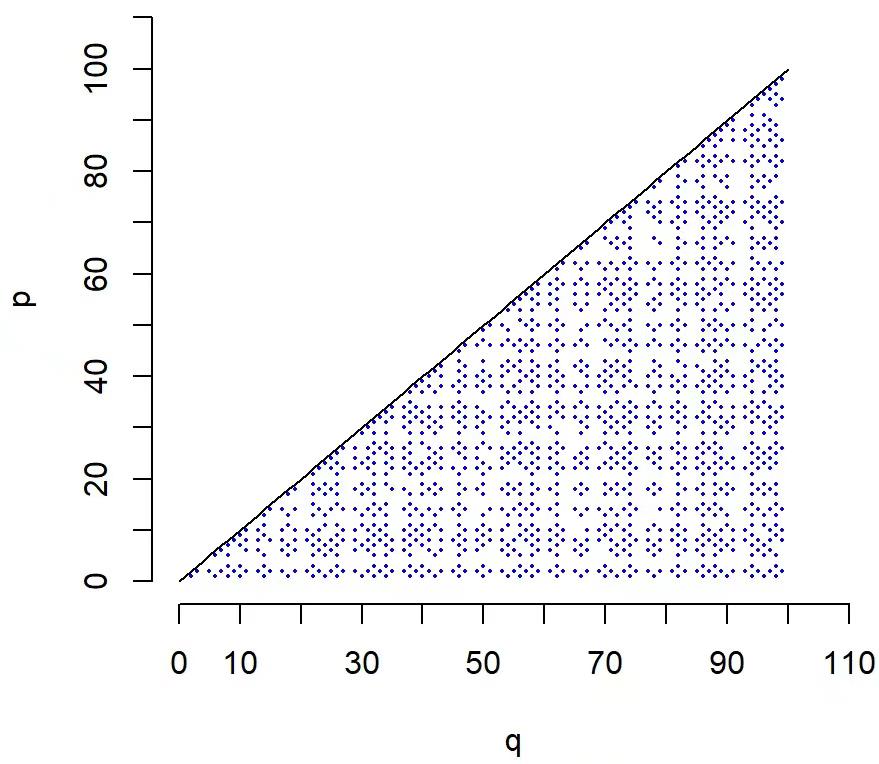}\\
  \caption{The graph of $\hat W\cap[1,N]^2$ with $N=100$.}\label{Fig:4}
\end{figure}
\begin{proof}
Note that for $n\in\N$, $ord_2(n)$ is odd if and only if $n=(2k-1)2^{2\ell-1}$ for some $k,\ell\in\N$.
  In view of  the definition of $\hat W$, we will count the number of pairs $(p,q)\in\hat W\cap[1, N]^2$  in the following way  (see Figure \ref{Fig:4}).
  First, by conditioned on $q=(2k-1)2^{2\ell-1}$  the number of   $p$'s satisfying  $(p,q)\in\hat W\cap[1, N]^2$ is $\phi((2k-1)2^{2\ell-1})$. Second, by conditioned on $p=(2k-1)2^{2\ell-1}$  the number of $q$'s satisfying   $(p,q)\in\hat W\cap[1, N]^2$ is given by $\varphi(N, (2k-1)2^{2\ell-1})-\phi((2k-1)2^{2\ell-1})$, where
  \[
  \varphi(N, (2k-1)2^{2\ell-1})=\#\set{1\le n\le N: n\textrm{ and }(2k-1)2^{2\ell-1}\textrm{ are coprime}}.
  \]
  For the range of $k$ and $\ell$, let $k_1$ be the largest $k\in\N$ such that $(2k-1)2^{2\ell-1}\le N$ for some $\ell\in\N$. Then
  \begin{equation}
 \label{eq:k1}
 k_1=k_1(N)=\flr{\frac{N+2}{4}}.
  \end{equation}
  Furthermore, for   $k\in[1, k_1]$ let $\ell_k$ be the largest $\ell\in\N$ such that $(2k-1)2^{2\ell-1}\le N$. Then
  \begin{equation}
 \label{eq:lk}
 \ell_k=\ell_k(N)=\flr{\log_4\Big(\frac{2N}{2k-1}\Big)}.
  \end{equation}
  Therefore,
  \begin{equation}
 \label{eq:Oct-1}
 \begin{split}
 \#(\hat W\cap[1,N]^2)&=\sum_{k=1}^{k_1}\sum_{\ell=1}^{\ell_k}\phi((2k-1)2^{2\ell-1})\\
 &\qquad+\sum_{k=1}^{k_1}\sum_{\ell=1}^{\ell_k}\Big(\varphi(N, (2k-1)2^{2\ell-1})-\phi((2k-1)2^{2\ell-1})\Big)\\
 &=\sum_{k=1}^{k_1}\sum_{\ell=1}^{\ell_k}\varphi(N,(2k-1)2^{2\ell-1}).
 \end{split}
  \end{equation}

 Observe that (cf.~\cite[Page 47, Exercise 9]{Apostol-1976})
 \begin{equation}
   \label{eq:Oct-2}
   \begin{split}
 \varphi(N, (2k-1)2^{2\ell-1})&=\sum_{n|(2k-1)2^{2\ell-1}}\mu(n)\flr{\frac{N}{n}}=\sum_{n|(2k-1)2^{2\ell-1}}\mu(n)\frac{N}{n}+\ep_{k,\ell}\\
 &=N\frac{\phi((2k-1)2^{2\ell-1})}{(2k-1)2^{2\ell-1}}+\ep_{k,\ell}=N\frac{\phi(2k-1)}{2(2k-1)}+\ep_{k,\ell},
   \end{split}
 \end{equation}
 where the last two equalities follow by (\ref{eq:Euler-function-1}) and (\ref{Euler's totient function}) respectively. Here the error term $\ep_{k,\ell}$ is bounded by
 \[
 |\ep_{k,\ell}|\le\sum_{n=1}^\f\mathbb I_{\set{n|(2k-1)2^{2\ell-1}}}=d((2k-1)2^{2\ell-1})=2\ell d(2k-1),
 \]
 where $d(m)$ denotes the number of all positive factors of $m$, and the last equality follows since $d(m)=\prod_{i=1}^n(c_i+1)$ if $m=\prod_{i=1}^n p_i^{c_i}$ (cf.~\cite[Theorem 273]{Hardy-Wright-2008}).
 Thus, by (\ref{eq:Oct-1}) and (\ref{eq:Oct-2}) it follows that
 \begin{equation}
   \label{eq:Oct-3}
   \begin{split}
 \lim_{N\to\f}\frac{\#(\hat W\cap[1,N]^2)}{N^2}&=\lim_{N\to\f}\frac{1}{N^2}\sum_{k=1}^{k_1}\sum_{\ell=1}^{\ell_k}\left(N\frac{\phi(2k-1)}{2(2k-1)}+\ep_{k,\ell}\right)\\
 &=\lim_{N\to\f}{\frac{1}{N}\sum_{k=1}^{k_1}\sum_{\ell=1}^{\ell_k}\frac{\phi(2k-1)}{2(2k-1)}},
   \end{split}
 \end{equation}
 where the last equality follows by (\ref{eq:k1}) and (\ref{eq:lk}) that
 \begin{align*}
   \left|\frac{1}{N^2}\sum_{k=1}^{k_1}\sum_{\ell=1}^{\ell_k}\ep_{k,\ell}\right|&\le\frac{1}{N^2}\sum_{k=1}^{k_1}\sum_{\ell=1}^{\ell_k}2\ell d(2k-1)=\frac{1}{N^2}\sum_{k=1}^{k_1}\ell_k(\ell_k+1)d(2k-1)\\
   &\le\frac{\log_4(2N)(\log_4(2N)+1)}{N^2}\sum_{k=1}^{\flr{\frac{N+2}{4}}}d(2k-1)\to 0\quad\textrm{as }N\to\f.
 \end{align*}
 Here the limit follows by \cite[Theorem 318]{Hardy-Wright-2008} that $\lim_{n\to\f}\frac{1}{n\ln n}\sum_{m=1}^n d(m)=1$.

 Note by (\ref{eq:k1}) that $k_1=\flr{\frac{N+2}{4}}$. In general, for $j\ge 1$ let
 \begin{equation}\label{eq:kj}
 k_j=k_j(N)=\flr{\frac{1}{2}+\frac{N}{4^j}}.
 \end{equation}
 Then $k_j$ decreases to zero as $j\to \f$.  In fact, for {$j\geq \flr{\log_4 2N}+1$} we have $k_j=0$. Note that for any $k\in(k_{j+1}, k_j]$ we have $\ell_k=j$. So, by (\ref{eq:Oct-3}) it follows that
 \begin{align*}
   \lim_{N\to\f}\frac{\#(\hat W\cap[1, N]^2)}{N^2}
   &=\lim_{N\to\f}\frac{1}{2N}\sum_{k=1}^{k_1}\ell_k\frac{\phi(2k-1)}{2k-1}\\
   &=\lim_{N\to\f}\frac{1}{2N}{\sum_{j=1}^{\flr{\log_4 2N}}\sum_{k=k_{j+1}+1}^{k_j}j\frac{\phi(2k-1)}{2k-1}}\\
   &=\lim_{N\to\f}\frac{1}{2N}{\sum_{j=1}^{\flr{\log_4 2N}}j\left(\sum_{k=1}^{k_j}\frac{\phi(2k-1)}{2k-1}-\sum_{k=1}^{k_{j+1}}\frac{\phi(2k-1)}{2k-1}\right).}
 \end{align*}
 Note by Lemma \ref{lem:sum-Euler-function-odd} that
 \[
 \sum_{k=1}^N\frac{\phi(2k-1)}{2k-1}=\frac{8}{\pi^2}N+o(N),
 \]where the little `$o$' stands for the higher order indefinite small.
 Therefore, by (\ref{eq:kj}) we obtain
 \begin{align*}
  \lim_{N\to\f} \frac{\#(\hat W\cap[1,N]^2)}{N^2}&=\lim_{N\to\f}\frac{1}{2N}\sum_{j=1}^{\flr{\log_4 2N}}j\left(\frac{8}{\pi^2}k_j-\frac{8}{\pi^2}k_{j+1}
   \right)\\
   &\qquad+\lim_{N\to\f}\frac{1}{2N}\sum_{j=1}^{\flr{\log_4 2N}}j\; o\big(k_j+k_{j+1}\big)\\
   &= \lim_{N\to\f}\frac{4}{\pi^2N}\sum_{j=1}^{\flr{\log_4 2N}}j\left(\flr{\frac{1}{2}+\frac{N}{4^j}}-\flr{\frac{1}{2}+\frac{N}{4^{j+1}}}\right)
   +\lim_{N\to\f}\frac{o(N)}{2N}\\
   &= \lim_{N\to\f}\frac{4}{\pi^2N}\sum_{j=1}^{\flr{\log_4 2N}}j\left( \frac{N}{4^j}- \frac{N}{4^{j+1}}\right)\\
   &=\frac{4}{\pi^2}\sum_{j=1}^\f\frac{3j}{4^{j+1}}=\frac{4}{3\pi^2}
 \end{align*}
  as desired.
\end{proof}

\begin{proof}
  [Proof of Theorem  \ref{th:density-exact-overlap}] By Proposition \ref{prop:density} we only need to consider the density of $W$. Note by (\ref{eq:W}) and (\ref{eq:hat-W}) that $W$ and $\hat W$ are disjoint, and
  \[
  W\cup\hat W=\set{(p, q)\in\N^2: p<q \textrm{ and }p, q\textrm{ are coprime}}.
  \]
  Then for large $N\in\N$ we have $\#(W\cup\hat W\cap[1,N]^2)=\sum_{n=2}^N\phi(n)$. So, by (\ref{eq:sum-Euler-function}) and Proposition \ref{prop:density} it follows that
  \begin{align*}
  \lim\limits_{N \to \infty}\frac{\#(W\cap [1,N]^{2})}{N^{2}}&=\lim\limits_{N \to \infty}\frac{\sum_{n=2}^{N}\phi (n)}{N^{2}}-\lim\limits_{N \to \infty}\frac{\#(\hat W\cap [1,N]^{2})}{N^{2}}\\
  &=\frac{3}{\pi^2}-\frac{4}{3\pi^2}=\frac{5}{3\pi^2},
  \end{align*}
completing the proof.
\end{proof}

\section*{Acknowledgements}
The first author was supported by Chongqing NSF: CQYC20220511052.

%

\end{document}